\documentclass[onefignum,onetabnum]{siamart190516}

\usepackage{lmodern}

\usepackage[utf8]{inputenc} 
\usepackage[T1]{fontenc}    
\usepackage{hyperref}       
\usepackage{url}            
\usepackage{booktabs}       
\usepackage{amsfonts}       
\usepackage{nicefrac}       
\usepackage{microtype}      
\usepackage{verbatim, amsfonts,enumerate}     
\usepackage{bm,array}
\usepackage{graphicx}
\usepackage{algorithm,algorithmic}
\usepackage{framed}
\usepackage{color,url}
\usepackage[normalem]{ulem}
\usepackage{float}
\usepackage{multirow}
\usepackage[font=normalsize]{caption}
\usepackage{hyperref}
\usepackage{subcaption}
\usepackage{cite}
\ifpdf
  \DeclareGraphicsExtensions{.eps,.pdf,.png,.jpg}
\else
  \DeclareGraphicsExtensions{.eps}
\fi
\graphicspath{{figure/}}

\newsiamremark{remark}{\bf Remark}
\newsiamremark{hypothesis}{Hypothesis}
\crefname{hypothesis}{Hypothesis}{Hypotheses}
\newsiamthm{claim}{Claim}
\newsiamthm{assumption}{Assumption}

\def \X {\mathbf{X}}
\def \Xt {{\tilde{\X}}}
\def \cX {\mathcal{X}}
\def \x {{\mathbf{x}}}
\def \A {{\mathbf{A}}}

\def \P {{\mathbf{P}}}
\def \F {{\mathcal{F}}}
\def \I {{\mathbf{I}}}

\def \a {{\mathbf{a}}}
\def \b {{\mathbf{b}}}
\def \e {{\mathbf{e}}}
\def \y {{\mathbf{y}}}
\def \yt {{\tilde{\y}}}

\def \c {{\mathbf{c}}}

\def \d {{\mathbf{d}}}

\def \v {{\mathbf{v}}}
\def \vt {{\tilde{\v}}}
\def \lam {{\bm{\lambda}}}

\def \Rn {\mathbb{R}}
\def \E {\mathcal{E}}
\def \O {\mathcal{O}}
\def \V {\mathcal{V}}
\def \N {\mathcal{N}}
\def \G {\mathcal{G}}
\def \Lam {\boldsymbol{\Lambda}}
\def \Lt {{\tilde{\Lam}}}

\def \lt {{\tilde{\lam}}}

\def \g {{\mathbf{g}}}

\def \h {\mathbf{h}}

\def \xt {{\tilde{\x}}}

\def \pr {{\text{prox}}}
\def \sX {{\mathsf{X}}}
\def \sY {{\mathsf{Y}}}
\def \mx {{\mu_{\sX}}}
\def \my {{\mu_{\sY}}}
\def \Dt {{\tilde{\Delta}}}
\def \bnu {{\boldsymbol{\nu}}}
\def \cD {{\mathcal{D}}}

\providecommand{\Et}[1]{\mathbb{E}_t\left[#1\right]}

\providecommand{\abs}[1]{\left|#1\right|}
\providecommand{\norm}[1]{\left\|#1\right\|}
\providecommand{\ip}[2]{\langle#1, #2\rangle}

\providecommand{\Ex}[1]{\mathbb{E}\left[#1\right]}
\headers{Primal-Dual Decentralized Stochastic Optimization}{K. Rajawat and C. Kumar}

\title{A Primal-Dual Framework for Decentralized Stochastic Optimization\thanks{Submitted to the editors DATE.
\funding{}}}

\author{ Ketan Rajawat\thanks{Department of Electrical Engineering, Indian Institute of Technology, Kanpur (India), 208016 
  (\email{ketan@iitk.ac.in} \url{}).}
\and Chirag Kumar\thanks{School of Electronics, Indian Institute of Information Technology Una, H.P.(India), 177220 
  (\email{chiragarora35@gmail.com}) .}}

\usepackage{amsopn}

\ifpdf
\hypersetup{
  pdftitle={A Primal-Dual Framework for Decentralized Stochastic Optimization},
  pdfauthor={K. Rajawat and C. Kumar}
}
\fi

\begin{document}

\maketitle

\begin{abstract}
  We consider the decentralized convex optimization problem, where multiple agents must cooperatively minimize a cumulative objective function , with each local function expressible as an empirical average of data-dependent losses. State-of-the-art approaches for decentralized optimization rely on gradient tracking, where consensus is enforced via a doubly stochastic mixing matrix. Construction of such mixing matrices is not straightforward and requires coordination even prior to the start of the optimization algorithm. This paper puts forth a primal-dual framework for decentralized stochastic optimization that obviates the need for such doubly stochastic matrices. Instead, dual variables are maintained to track the disagreement between neighbors. The proposed framework is flexible and is used to develop decentralized variants of SAGA, L-SVRG, SVRG++, and SEGA algorithms. Using a unified proof, we establish that the oracle complexity of these decentralized variants is $O(1/\epsilon)$, matching the complexity bounds obtained for the centralized variants. Additionally, we also present a decentralized primal-dual accelerated SVRG algorithm achieving $O(1/\sqrt{\epsilon})$ oracle complexity, again matching the bound for the centralized accelerated SVRG. Numerical tests on the algorithms establish their superior performance as compared to the variance-reduced gradient tracking algorithms.
\end{abstract}



\section{Introduction}
\label{intro}
Modern consumer devices generate massive amounts of data and successful deployment of large-scale machine learning algorithms necessitates addressing the challenges of scalability and privacy. Decentralized optimization offers a viable approach towards the design of such learning platforms. Within the decentralized learning rubric, each node retains the user data, while only iterates are exchanged between neighbors. Decentralized optimization algorithms exploit the availability of on-board device processors and therefore work in a peer-to-peer fashion without any central coordinator or parameter server. 

In this work, we consider $V$ nodes that seek to cooperatively solve the following network optimization problem:
\begin{align}\label{introprob}
\min_{\x \in \cX_i} \sum_{i=1}^V (f_i(\x) + h_i(\x)) 
\end{align}
where the data-dependent convex loss function $f_i$ is local to a node $i$ and takes the finite-sum form $f_i(\x) := \frac{1}{K_i}\sum_{k=1}^{K_i} f_i(\x,\xi_i^k)$, with $\{\xi_i^k\}_{k=1}^{K_i}$ denoting data points, which could be images, videos, or browsing history, stored at node $i$. The regularization function $h_i$ and the constraint set $\cX_i$ are also convex. The network topology is denoted by a connected graph $\G$ and the nodes may only communicate over the edges of $\G$. 

\begin{table*}	[thb]
\caption{Classes of decentralized algorithms: references within a cell are not exhaustive but only representative of the corresponding class. A cell mentions "yes" if at least one algorithm in the class exhibits the corresponding property. }\label{cat}
{\small
	\centering
	{\renewcommand{\arraystretch}{1.1}
		\begin{tabular}{|c|c|c|c|c|}
			\hline
			References  & Stochastic & \multirow{2}{*}{\begin{tabular}{c} Projected or \\
					proximal variant\end{tabular}} & \multirow{2}{*}{\begin{tabular}{c} Doubly stochastic \\
					mixing matrix\end{tabular}} & \multirow{2}{*}{\begin{tabular}{c} Variance \\
					reduced\end{tabular}}\\
			& & & & \\
			\hline
			\multirow{2}{*}{\begin{tabular}{c}
					\cite{tian2018achieving,qu2017harnessing,in2019distributed}\\
					\cite{qu2019accelerated,alghunaim2019linearly, xu2019accelerated}
			\end{tabular}}   & & & &  \\
			& No &  Yes & No &  No\\ 
			\hline
			\multirow{2}{*}{\begin{tabular}{c}
					\cite{xin2019distributed,pu2018push,pu2018distributed} \\
					\cite{li2019communication,xin2019variance,li2019s}
			\end{tabular}}   & & & &  \\
			& Yes &  No & Yes &  Yes \\
			\hline
			\cite{lan2018communication} & Yes & Yes & No & No \\
			\hline
			\cite{wang2019edge} & Yes & No & No & Yes \\
			\hline
			Proposed & Yes & Yes & No & Yes \\
			\hline
	\end{tabular}}}
\end{table*}
The distributed (sub-)gradient descent (DGD) and its variants were first proposed to solve the deterministic version of \eqref{introprob}; see \cite{nedic2018distributed} and references therein. The deterministic version of \eqref{introprob} has been well-studied, and a number of approaches achieving optimal rates exist; see \cite{tian2018achieving,qu2017harnessing,in2019distributed,qu2019accelerated,alghunaim2019linearly, xu2019accelerated} and others. Subsequent years saw the development of stochastic variants of DGD, with the state-of-the-art decentralized algorithms relying on the idea of gradient tracking \cite{di2015distributed,di2016next,xu2015augmented,xu2017convergence,xin2019distributed,pu2018push,pu2018distributed,li2019communication,xin2019variance,li2019s,sun2019convergence}. While most of these algorithms have only been analyzed for strongly convex objectives, they achieve optimal rates matching those obtained in centralized settings. However, stochastic gradient tracking algorithms suffer from two critical issues that render them unattractive for solving \eqref{introprob} in a scalable fashion: (a) projected or proximal variants of stochastic gradient tracking do not exist and cannot be analyzed using existing proof techniques that rely on unrolling the iterations; and (b) one or more doubly stochastic mixing matrices are required to be used for information fusion between nodes. As a result, the oracle complexity bounds obtained for gradient tracking algorithms always depend on the spectral properties of these mixing matrices. In practice, the mixing matrices are often designed in an ad hoc fashion and necessitate network-wide coordination even before the start of the algorithm. Likewise, the requirements that $h_i = 0$ and $\cX = \Rn^n$ severely limit the class of optimization problems that can be solved using gradient tracking methods.

Primal-dual class of algorithms are capable of addressing the aforementioned challenges to a certain extent. Here, the consensus constraints are written explicitly and Lagrange multipliers are utilized to ensure network-wide consensus. A primal-dual algorithm to solve \eqref{introprob} was proposed in \cite{lan2018communication} where distributed communication sliding (DCS) was utilized to reduce the inter-node communications. However, the stochastic version of DCS proposed in \cite{lan2018communication} was primarily meant for online settings and enhancements such as variance reduction were not applied, thereby resulting in an oracle complexity of $\O(\frac{1}{\epsilon^2})$ for general convex functions. An edge-based primal-dual stochastic average gradient method is proposed in \cite{wang2019edge}, and is somewhat similar to the algorithms considered here, but is analyzed only for strongly convex problems. In contrast, the present work focuses on the general convex and smooth functions, for which there exist no decentralized stochastic algorithms achieving optimal rates. 

Beyond these, there also exist decentralized algorithms that require oracles other than the stochastic gradient oracle. For instance, the updates in \cite{hendrikx2019accelerated} require the proximal operator to be calculated for each component $f_i(\cdot,\xi_i^k)$. The works in \cite{dvinskikh2019primal,gorbunov2019optimal} can also solve decentralized problems but require access to stochastic dual oracles. However, such oracles may not easily be available for many machine learning problems. Table \ref{cat} summarizes the different classes of decentralized algorithms for convex optimization problems.

This work puts forth a primal-dual framework for the construction of decentralized stochastic variance-reduced algorithms for solving \eqref{introprob}. The proposed framework is accompanied by a unified analysis technique that allows us to construct decentralized versions of several variance-reduction techniques. In the present work, we propose decentralized variants of the SVRG++ \cite{allen2016improved}, loopless SVRG \cite{kovalev2019don}, SAGA \cite{defazio2014saga}, and SEGA \cite{gorbunov2019unified}, though the framework allows construction of other variants as well. We show that similar to the centralized versions, these decentralized algorithms also achieve an oracle complexity of $\O(\frac{1}{\epsilon})$ for general convex functions. 

Additionally, and for the first time in the context of stochastic decentralized optimization, we propose an accelerated decentralized variance reduced primal-dual algorithm with oracle complexity $\O(\frac{1}{\sqrt{\epsilon}})$. The result is remarkable as it improves upon the existing decentralized algorithms working with stochastic gradient oracles. Unlike gradient tracking algorithms, the proposed algorithms can handle arbitrary convex regularizers and constraints, and do not require doubly stochastic mixing matrices.  

Table \ref{oracle} summarizes the various oracle and communication complexity results obtained in this paper. For the communication complexity, we count the total number of network-wide communication rounds. All the results are novel and have not been reported earlier. In the table, $L$ denotes the smoothness parameter while $K$ denotes the number of data points per node; detailed results can be found in Sec. \ref{variants} and \ref{accsec}.

This paper is organized as follows: Sec. \ref{pdalgos} develops the proposed primal-dual algorithm in its generic form. Sec. \ref{variants} provides the stochastic gradient, SVR, SAGA, and LSVR variants of the proposed algorithm and the corresponding oracle complexity results. Sec. \ref{accsec} develops the accelerated variant of the decentralized SVRG and provides the corresponding oracle complexity result. 

\emph{Notation:} Vectors are denoted by small bold letters, matrices by capital bold letters, and sets by calligraphic capital letters. All norms are $\ell_2$ norms. 


\begin{table}[]
	\caption{Oracle complexities of the proposed algorithms for general convex and $L$-smooth losses. The $\O(\cdot)$ notation subsumes initialization-dependent quantities. The network has $V$ nodes and $E$ edges, and each node has $K$ data points.}\label{oracle}
	\centering
	{\renewcommand{\arraystretch}{1.8}
	\begin{tabular}{|m{3.5cm}|m{4cm}|m{3cm}|}
		\hline
		{\renewcommand{\arraystretch}{1.2}\begin{tabular}[c]{@{}l@{}}Decentralized \\ Algorithm\end{tabular}} & 
		\begin{tabular}[c]{@{}l@{}}Oracle Complexity\\ \end{tabular} & {\renewcommand{\arraystretch}{1.2}\begin{tabular}[c]{@{}l@{}}Communication \\ Complexity\end{tabular}} \\ 
		\hline
		SAGA-PD                                                            & \multicolumn{2}{c|}{\multirow{2}{*}{$\O\left(\frac{E+V\min\{L+K,\sqrt{LK}\}}{\epsilon}\right)$}}                                                                                        \\ \cline{1-1}
		LSVR-PD & \multicolumn{2}{c|}{} \\ 
		\hline
		SVR-PD & $\O(\frac{E+VL}{\epsilon} + K\log(\frac{E+VL}{\epsilon}))$ &  $\O(\frac{E+VL}{\epsilon})$ \\ 
		\hline
		\begin{tabular}[c]{@{}l@{}}ASVR-PD\\ \end{tabular} & $\O\left(\frac{\sqrt{K(E+V(K+L))}}{\sqrt{\epsilon}}\right)$   &  $\O(\frac{E+V(K+L)}{\epsilon})$  \\ 
		\hline
	\end{tabular}}
\end{table}

\def \colb {\color{blue}}
\section{Decentralized Stochastic Primal-Dual Algorithm}\label{pdalgos}
Consider a multi-agent network represented by an undirected graph $\G = (\V, \E)$, where $\V$ denotes the set of nodes and $\E=\{(i,j) | i,j \in \V\}$ denotes the set of edges. There are a total of $V$ nodes and $E$ edges. The set of neighboring nodes of $i$ are denoted by $\N_i:=\{j \in \V \mid (i,j) \in \E\}$ and its cardinality is given by $\abs{\N_i}$. A node $i$ is associated with a proper closed convex and smooth function $f_i:\Rn^n\rightarrow\Rn$. The network seeks to cooperatively solve the  optimization problem $\min_{\x} \sum_{i=1}^V (f_i(\x) + h_i(\x))$ where $\x \in \Rn^n$ is the optimization variable and $h_i:\Rn^n \rightarrow \Rn$ are proper closed convex but possibly non-smooth regularization functions. Constraints of the form $\x\in\cX_i$ where $\cX_i$ is closed and convex, if any, are included as indicator functions $1_{\x\in\cX_i}$ within $h_i$. 
As discussed earlier, the loss function takes the form
\begin{align}
f_i(\x) = \frac{1}{K_i}\sum_{k=1}^{K_i} f_i(\x,\xi_i^k)\label{finitesum}
\end{align}
for all $i\in\V$. Throughout the paper, we will work with the following key assumption.

\begin{assumption}\label{a1}
	The functions $f_i(\cdot, \xi_i)$ are proper, closed, convex, and $L_i$-smooth. The functions $h_i$ are proper, closed, and convex. 
\end{assumption}

\subsection{Generic Decentralized Primal-Dual Algorithm}
We begin with describing the generic version of the proposed primal-dual (PD) algorithm. The specific choice of the stochastic gradients corresponding to different variance-reduced versions and their effect on the oracle complexity will be studied in the subsequent sections. Associate local variable $\x_i\in\Rn^n$ for each node $i\in\V$, all collected into the matrix $\X \in \Rn^{n \times V}$. The overall problem can now be written as
\begin{align}
\X^\star& = \arg\min_{\X} \sum_{i=1}^V [f_i(\x_i) + h_i(\x_i)]  \hspace{1cm}\text{  s. t. }  \  \x_i = \x_j \hspace{4mm} \forall(i,j) \in \E \tag{\mbox{$\mathcal{P}$}}.
\label{P1}
\end{align}
The domain of the objective function is denoted by $\cD:=\{\text{dom}f_1 \cap \text{dom}h_1\} \times \ldots \times \{\text{dom}f_V \cap \text{dom}h_V\}$ and is assumed to have a non-empty relative interior. 
The standard augmented Lagrangian for \eqref{P1}, such as that used in the context of ADMM algorithms, is given by 
\begin{align}
\sum_{i\in\V} (f_i(\x_i) + h_i(\x_i)) &+ \sum_{(i,j)\in\E}\ip{\lam_{ij}}{\x_i-\x_j} + \rho\sum_{(i,j)\in\E}\norm{\x_i-\x_j}^2.\label{aug}
\end{align}
Note that since each constraint is written only once, only one of $\lam_{ij}$ or $\lam_{ji}$ appears in \eqref{aug}. To make the algorithm distributed, we instead consider two dual variables $\lam_{ij}$ and $\lam_{ji}$ per edge. Under the dual constraint $\lam_{ij}+\lam_{ji} = 0$, the augmented Lagrangian can be written as 
\begin{align}
L(\X,\Lam) = \sum_{i\in\V} (f_i(\x_i) +& h_i(\x_i)) + \sum_{i\in\V}\sum_{j\in\N_i}\ip{\lam_{ji}}{\x_i} + \frac{\rho}{2}\sum_{i\in\V}\sum_{j\in\N_i}\norm{\x_i-\x_j}^2.\label{augl}
\end{align}
where the matrix $\Lam \in \Lambda \subset  \Rn^{n \times 2\abs{\E}}$ collects all the dual variables with $\Lambda =  \{\Lam | \lam_{ij}+\lam_{ji} = 0\}$. The domain of $L(\X,\Lam)$ is given by  dom$(L(\X,\Lam)) = \cD \times \Lambda$. The following lemma, whose proof follows from standard duality theory arguments, and is provided in Appendix \ref{exactgrad},  
lists the optimality conditions of \eqref{P1} and relates them with the Lagrangian $L(\X,\Lam)$. 
\begin{lemma}\label{opt}
Given  $\X^\star$, there exists $\Lam^\star$ such that $(\X^\star,\Lam^\star)$ is a saddle point of $L(\X,\Lam)$ and satisfies the following conditions:
	\begin{subequations}\label{optcon}
		\begin{align}
		\x_i^\star - \x_j^\star &= 0 \label{pfeas}\\
		\lam_{ij}^\star + \lam_{ji}^\star &= 0 \label{dfeas}\\
		\v_i^\star + \nabla f_i(\x_i^\star) + \sum_{j\in\N_i}\lam_{ji}^\star &= 0. \label{fo}
		\end{align}
	\end{subequations}
for $\v_i^\star \in \partial h_i(\x_i^\star)$.
\end{lemma}
The proximal operator is defined as $\pr_{\eta h}(\y) = \arg\min_{\x} h(\x) + \frac{1}{2\eta}\norm{\y-\x}^2$. For constants $\eta_i > 0$ and $\rho >0$,  Lemma \ref{opt} implies that 
\begin{subequations}\label{optcon2}
\begin{align}
\x_i^\star &= \pr_{\eta_i h_i}\Big(\x_i^\star - \eta_i\nabla f_i(\x_i^\star) + \eta_i\sum_{j\in\N_i}\lam_{ij}^\star\Big) \\
\lam_{ij}^\star &= -\lam_{ji}^\star + \rho(\x_j^\star - \x_i^\star).
\end{align}\end{subequations}
The optimality conditions in \eqref{optcon2} suggest the following proximal stochastic saddle point algorithm:
\begin{align}
\x_i^{t+1} &= \pr_{\eta_i h_i}\Big(\x_i^t - \eta_i\g_i^t + \eta_i\sum_{j\in\N_i}\lam_{ij}^{t+1}\Big) \label{primup1}\\
\lam_{ij}^{t+1} &= -\lam_{ji}^t + \rho(\x_j^t - \x_i^{t+1}). \label{dualup1}
\end{align}
where $\g_i^t$ is a stochastic approximation of $\nabla f_i(\x_i^t)$. From \eqref{optcon2}, it can be seen that when exact gradient is used (i.e. $\g_i^t = \nabla f_i(\x_i^t)$), the primal-dual optimum $(\X^\star,\Lam^\star)$ is a limit point of \eqref{primup1}-\eqref{dualup1}. Further, as in the proof of Lemma \ref{opt}, it can be established that any limit point of \eqref{primup1}-\eqref{dualup1} satisfies \eqref{optcon} in the exact gradient case; see Appendix \ref{exactgrad}.

 An implementable form for the updates in \eqref{primup1}-\eqref{dualup1} is obtained by substituting $\lam_{ij}^{t+1}$ from \eqref{dualup1} into \eqref{primup1} and simplifying, which yields:
\begin{subequations}
	\begin{align}
	\x_i^{t+1} &= \pr_{\gamma_i h_i}\Big(\frac{\gamma_i}{\eta_i}[\x_i^t-\eta_i\g_i^t] + \gamma_i\sum_{j\in\N_i}[\rho\x_j^t-\lam_{ji}^t] \Big) \label{primup}\\
	\lam_{ij}^{t+1} &= -\lam_{ji}^t + \rho(\x_j^t - \x_i^{t+1})\label{dualup}
	\end{align}
\end{subequations}
for all $j\in \N_i$ and $i\in\V$, where $\gamma_i = \big(\frac{1}{\eta_i} + \rho\abs{\N_i}\big)^{-1}$. Here, $\abs{\N_i}$ is the number of neighbors of node $i$. The updates can be carried out in parallel at each node $i$, given $\{(\lam_{ji}^t, \x_j^t)\}_{j\in\N_i}$. 
The generic decentralized primal-dual algorithm is summarized in Algorithm \ref{generic}. Observe that Algorithm \ref{generic} requires node $i$ to maintain $\{\lam_{ij}^{t}\}_{j\in\N_i}$ and obtain $\{(\x_j^t, \lam_{ji}^t)\}_{j\in\N_i}$ from its neighbors. After the updates, the node broadcasts $\{(\x_i^{t+1},\lam_{ij}^{t+1})\}_{j\in\N_i}$ to all its neighbors. 
\begin{algorithm}
	\caption{Decentralized Primal-Dual}\label{generic}
	\begin{algorithmic}
	\STATE Initialize $\{\x^1_i\}_{i\in\V}$,  $\{\lam^1_{ij}\}_{i\in\V}$
	\FOR { $\textit{t} = 1,2,... $}
		\FOR{ $i \hspace{2mm} \in \hspace{2mm} \V$}
				\STATE Approximate $\g_i^t \approx \nabla f_i(\x_i^t)$
				\STATE Update $\x_i^{t+1}$ using \eqref{primup}
				\STATE Update $\lam_{ij}^{t+1}$ using \eqref{dualup}	 	\ENDFOR
\ENDFOR
\end{algorithmic}
\end{algorithm}

When using the exact gradient, i.e., $\g_i^t = \nabla f_i(\x_i^t)$, the decentralized PD algorithm is related to the quadratic approximation primal-dual method of multipliers (PDMM) algorithm proposed in \cite{qapdmm}; see Appendix \ref{PDMM} for details. However, the quadratic approximation PDMM has only been studied in the deterministic setting for the case $h_i = 0$. The analysis presented here is considerably more general and applicable to large-scale settings. Next we state some preliminaries and develop a general bound on the optimality gap for Algorithm \ref{generic}, which will be used in the next section to study the different variants of the proposed algorithm.

\subsection{Oracle Complexity Analysis: Preliminaries}
We begin with deriving some of the basic inequalities that will be used throughout the proofs. The Peter-Paul inequality states that for any $\a$, $\b \in \Rn^n$, and scalar $c > 0$, $2\ip{\a}{\b} \leq c \norm{\a}^2 + \frac{1}{c}\norm{\b}^2$ which, for $c = 1$, translates to $\norm{\a+\b} \leq 2\norm{\a}^2 + 2\norm{\b}^2$. For two random variables $\mathsf{X}$ and $\mathsf{Y}$, since variance is non-negative, we have that
\begin{align}
\Ex{\norm{\sX - \Ex{\sX} + \sY}^2} &\leq 2\Ex{\norm{\sX-\Ex{\sX}}^2} + 2\Ex{\norm{\sY}^2} \nonumber\\
&\leq 2\Ex{\norm{\sX}^2} + 2\Ex{\norm{\sY}^2}.\label{expineq}
\end{align}

In order to study the performance of the proposed algorithms, we make use of the Bregman divergences between two points $\x$ and $\y\in\cX$ evaluated along the convex functions $f_i$ and $h_i$, 
\begin{align}
D_{f_i}(\x,\y)&:=f_i(\x) - f_i(\y) - \ip{\nabla f_i(\y)}{\x-\y}\label{brg}\\
D_{h_i}(\x,\y)&:=h_i(\x) - h_i(\y) - \ip{\v_i(\y)}{\x-\y}
\end{align}
where $\v_i(\y) \in\partial h_i(\y)$. Recall that Bregman divergences are non-negative but need not be symmetric. Note that here we consider a generalized definition of Bregman divergence for any proper, closed, convex functions $f_i$ and $h_i$. If $f_i$ is also strictly convex, the condition $D_{f_i}(\x,\y) = 0$ would be equivalent to $\x = \y$. In general however, we only have that 
\begin{align}
\ip{\x-\y}{\nabla f_i(\x) - \nabla f_i(\y)} \geq D_{f_i}(\x,\y) \geq \frac{1}{2L_i}\norm{\nabla f_i(\x) - \nabla f_i(\y)}^2. \label{dfineq}
\end{align} 
Further, since $f_i(\cdot,\xi_i^k)$ is smooth from Assumption \ref{a1}, it follows from \eqref{dfineq}
that
\begin{align}
\frac{1}{K_i}\sum_{k=1}^{K_i}\norm{\nabla f_i(\x_i^t,\xi_i^k) - \nabla f_i(\x_i^\star,\xi_i^k)}^2\leq 2L_iD_{f_i}(\x_i^t,\x_i^\star). \label{dfbound}
\end{align}	
where we have used the definition of $f_i$ in \eqref{finitesum}. The following non-negative function will be useful in the analysis of the proposed algorithm
\begin{align}
\Psi^t :=  \frac{1}{2\rho}\sum_{i\in\V}\sum_{j\in\N_i}\norm{\rho(\x_i^t-\x_i^\star) - (\lam_{ij}^t-\lam_{ij}^\star)}^2. \label{bk}
\end{align}
We assume that the algorithm is initialized such that 
\begin{align}
	\max_{i,j}\{\norm{\x_i^1-\x_i^\star}^2, \norm{\lam_{ij}^1 -\lam_{ij}^\star}^2, D_{f_i+h_i}(\x_i^1,\x_i^\star)\} &\leq W_1 \label{w1}
\end{align}
so that $\Psi^1 \leq 2E(\rho + \tfrac{1}{\rho})W_1$ where recall that $E = \abs{\E}$ and $W_1$ is some constant. The complexity of all algorithms will be specified in terms of the number of oracle calls required to reach a pre-specified accuracy. The optimality gap of the different algorithmic variants will be characterized using the cumulative Bregman divergence given by 
\begin{align}\label{bregman}
\Delta^t = \sum_{i\in\V} \Ex{D_{f_i+h_i}(\x_i^t,\x^\star)}
\end{align}
where the expectation is taken with respect to random gradients $\{\g_i^1, \ldots, \g_i^t\}_{i\in\V}$. In particular, we will obtain bounds on the number of oracle calls, i.e., evaluations of $\{\nabla f_i(\x_i^t,\xi_i^k)\}_{i\in\V}$, required to ensure that $\min_{1\leq t \leq T} \Delta^t \leq \epsilon$. A small value of $\Delta^t$ implies that $\x_i^t$ and $\x_i^\star$ are close along $f_i+h_i$, for all $i$. While the bounds are developed for general case, we will often drop the node index and replace $L_i$ with $L$ and $K_i$ with $K$ to get intuition and express the results simply. Communication complexity is given by the total number of communication rounds over the entire horizon.

\subsection{Analysis of Algorithm \ref{generic}}
The oracle complexity analysis of the non-accelerated algorithms will follow a uniform template similar to that in \cite{gorbunov2019unified}, except that we consider general (non-strongly) convex functions here. We begin with stating a basic result for the proposed class of algorithms. 

\begin{lemma}\label{lema1}
	The function difference $\Psi^{t+1}-\Psi^t$ satisfies
	\begin{align}
	\Psi^{t+1}-\Psi^t  &= -2\sum_{i\in\V,j\in\N_i}\ip{\x_i^{t+1} - \x_i^\star}{\lam_{ij}^{t+1}-\lam_{ij}^\star}\label{lyap1}\\
	&\leq -\frac{2}{\eta_i}\sum_{i\in\V} \ip{\x_i^{t+1} - \x_i^\star}{\x_i^{t+1}-\x_i^t} - 2\sum_{i\in\V}\ip{\x_i^{t+1} - \x_i^\star}{\g_i^t-\nabla f_i(\x_i^\star)}\label{lyap2}\\
	& - 2\sum_{i\in\V}D_{h_i}(\x_i^{t+1},\x_i^\star)\nonumber
	\end{align}
\end{lemma}
\begin{proof}[Proof of Lemma \ref{lema1}]	
Subtracting $\lam_{ij}^\star + \rho\x_i^\star$ from both the sides of \eqref{dualup}, using the facts that $\lam_{ij}^\star + \lam_{ji}^\star = 0$ and  $\x_i^\star = \x_j^\star$, and rearranging, we obtain
\begin{align}
\lam_{ij}^{t+1} - \lam_{ij}^\star +  \rho(\x_i^{t+1}-\x_i^\star) = -(\lam_{ji}^t - \lam_{ji}^\star) + \rho(\x_j^t - \x_j^\star).\label{dualup-alt}
\end{align}                         
Using \eqref{dualup-alt}, the $(i,j)$-th summand on the right of \eqref{lyap2} can be written as
\begin{align}
&\frac{1}{2\rho}\norm{\rho(\x_i^{t+1}-\x_i^\star) - (\lam_{ij}^{t+1}-\lam_{ij}^\star)}^2 - \frac{1}{2\rho}\norm{\rho(\x_i^t-\x_i^\star) - (\lam_{ij}^t-\lam_{ij}^\star)}^2 \nonumber\\
&= \frac{1}{2\rho}\norm{\rho(\x_i^{t+1}-\x_i^\star) - (\lam_{ij}^{t+1}-\lam_{ij}^\star)}^2 - \frac{1}{2\rho}\norm{\rho(\x_j^{t+1} - \x_j^\star) + \lam_{ji}^{t+1} - \lam_{ji}^\star}^2 \nonumber\\
& = \frac{\rho}{2}\norm{\x_i^{t+1} - \x_i^\star}^2 - \frac{\rho}{2}\norm{\x_j^{t+1} - \x_j^\star}^2 + \frac{1}{2\rho}\norm{\lam_{ij}^{t+1}-\lam_{ij}^\star}^2 - \frac{1}{2\rho}\norm{\lam_{ji}^{t+1} - \lam_{ji}^\star}^2 \nonumber\\
&  - \ip{\x_i^{t+1} - \x_i^\star}{\lam_{ij}^{t+1}-\lam_{ij}^\star} - \ip{\x_j^{t+1} - \x_j^\star}{\lam_{ji}^{t+1} - \lam_{ji}^\star}.\nonumber
\end{align}
Summing both sides over $i\in\V$ and $j\in \N_i$, observe that the first four terms cancel, and we obtain
\begin{align}
\Psi^{t+1}-\Psi^t  &= -2\sum_{i\in\V,j\in\N_i}\ip{\x_i^{t+1} - \x_i^\star}{\lam_{ij}^{t+1}-\lam_{ij}^\star}
\end{align}
as required in \eqref{lyap1}. 
Next, the update equation in \eqref{primup1} and the optimality conditions in \eqref{optcon} imply that
\begin{align}
\x_i^{t+1} &= \x_i^t - \eta_i\g_i^t + \eta_i\sum_{j\in\N_i}\lam_{ij}^{t+1} - \eta_i\v_i^{t+1} \\
\Rightarrow \sum_{j\in\N_i}(\lam_{ij}^{t+1} - \lam_{ij}^\star) &= \tfrac{1}{\eta_i}(\x_i^{t+1} - \x_i^t) + \g_i^t - \nabla f_i(\x_i^\star) + \v_i^{t+1} - \v_i^\star 
\end{align}
where $\v_i^{t+1} \in \partial h_i(\x_i^{t+1})$. Upon substituting in \eqref{lyap1}, we obtain
\begin{align}
\Psi^{t+1}-\Psi^t &= -\frac{2}{\eta_i}\sum_{i\in\V} \ip{\x_i^{t+1} - \x_i^\star}{\x_i^{t+1}-\x_i^t} - 2\sum_{i\in\V}\ip{\x_i^{t+1} - \x_i^\star}{\g_i^t-\nabla f_i(\x_i^\star)}\\ &- 2\sum_{i\in\V}\ip{\x_i^{t+1} - \x_i^\star}{\v_i^{t+1}-\v_i^\star}.\nonumber \label{last1}
\end{align}
 Since $h_i$ is convex, each summand in the last term can be bounded as
\begin{align}
\ip{\x_i^{t+1} - \x_i^\star}{\v_i^{t+1}-\v_i^\star} \geq h_i(\x_i^{t+1}) - h_i(\x_i^\star) - \ip{\v_i^\star}{\x_i^{t+1}-\x_i^\star} = D_{h_i}(\x_i^{t+1},\x_i^\star).\nonumber
\end{align} 
Substituting, we obtain the desired result in \eqref{lyap2}. 
\end{proof}

We make the following assumption regarding the first and second moments of $\g_i^t$.  Let $\F_t$ denote the sigma algebra formed by the random variables $\{\g_i^1, \ldots \g_i^{t-1}\}_{i\in\V}$ and for any random variable $\mathsf{X}$, let $\Et{\mathsf{X}} := \Ex{\mathsf{X}\mid \F_t}$.  
\begin{assumption}\label{a2}
	The stochastic gradient $\g_i^t$ is conditionally unbiased given $\F_t$, i.e., 
	\begin{align}
	\Et{\g_i^t} = \nabla f_i(\x_i^t) \label{unbiased}
	\end{align}
	and there exist non-negative constants $A_i$, $B_i$, $C_i$, $D_1$, $D_2$, $\sigma_i^2$, $\varrho_i > 0$, and a random sequence $\phi_i^t \geq 0$ such that the following inequalities hold:
	\begin{align}
	\Et{\norm{\g_i^t -\nabla f_i(\x_i^\star)}^2} &\leq 2A_iD_{f_i}(\x_i^t,\x^\star_i) + B_i \phi_i^t + 2\sigma_i^2\label{a11}\\
	\Et{\phi_i^{t+1}} &\leq (1-\varrho_i)\phi_i^t + 2C_iD_{f_i}(\x_i^t,\x^\star_i)	\label{a22}
	\end{align}
\end{assumption}
The form of Assumption \ref{a2} is motivated from \cite{gorbunov2019unified} where it was used to analyze different variance-reduced algorithms under the strongly-convex setting. In the present case also, the assumption allows us to study the different SGD variants in a unified manner. Different from \cite{gorbunov2019unified}, we consider general convex functions. The following lemma provides the convergence rate of Algorithm \ref{generic} and will form the basis of our analysis of various stochastic algorithms. 

\begin{lemma}\label{lema2}
	Under Assumptions \ref{a1} and \ref{a2}, it holds that
	\begin{align} \label{mainbound}
	\frac{1}{T}\sum_{t=1}^T \Delta^t &\leq \left(2E\left(\rho+\frac{1}{\rho}\right) + \sum_{i\in\V} \frac{1}{\eta_i} + 2V\right)\frac{W_1}{T} + \sum_{i\in\V}\frac{B_i\eta_i\phi_i^1}{\varrho_i T} + 2\sum_{i\in\V} \sigma_i^2\eta_i \\
	&= \frac{W_1G_\rho}{T} + \sum_{i\in \V}\left[\frac{W_1}{\eta_i T} + \frac{B_i\eta_i\phi_i^1}{\varrho_i T} + \sigma_i^2\eta_i\right]
	\end{align}
	where $\eta_i \leq (2A_i + \frac{2B_iC_i}{\varrho_i})^{-1}$ for all $i\in\V$, $W_1$ defined as in \eqref{w1}, and $G_\rho = 2E(\rho+1/\rho)+2V$. 
\end{lemma}

\begin{proof}
Taking conditional expectation in \eqref{lyap2}, we obtain
\begin{align}
\Et{\Psi^{t+1} -\Psi^t}\leq& -\frac{2}{\eta_i}\sum_{i\in\V} \Et{\ip{\x_i^{t+1} - \x_i^\star}{\x_i^{t+1}-\x_i^t}}\\  &-2\sum_{i\in\V}\ip{\x_i^t-\x_i^\star}{\nabla f(\x_i^t) -\nabla f_i(\x_i^\star)}\nonumber\\
&- 2\sum_{i\in\V}\Et{\ip{\x_i^{t+1}-\x_i^t}{\g_i^t-\nabla f_i(\x_i^\star)}}\nonumber\\  &- 2\sum_{i\in\V}\Et{D_{h_i}(\x_i^{t+1},\x_i^\star)}\nonumber\\
& =: \sum_{i\in\V} ( I_1^i + I_2^i + I_3^i + I_4^i) 
\end{align}
where $I_1^i$, $I_2^i$, $I_3^i$, and $I_4^i$ denote the $i$-th summand in the first, second, third, and fourth terms, respectively. For the first term observe that
\begin{align}
\eta_i I_1^i &= -2\Et{\ip{\x_i^{t+1}-\x_i^\star}{\x_i^{t+1}-\x_i^t}} \\
&= \norm{\x_i^t-\x_i^\star}^2 - \Et{\norm{\x_i^{t+1}-\x_i^\star}^2}  - \Et{\norm{\x_i^{t+1}-\x_i^t}^2}. \label{i1}
\end{align}
Since $f_i$ is convex, the second term can be bounded from \eqref{dfineq} as follows:
\begin{align}
I_2^i &= -2\ip{\x_i^t-\x_i^\star}{\nabla f_i(\x_i^t)-\nabla f_i(\x_i^\star)} \\
&\leq -2D_{f_i}(\x_i^t,\x_i^\star).\label{i2}
\end{align}
We use the Peter-Paul inequality with parameter $\eta_i$ to bound $I_3^i$ as
\begin{align}
I_3^i &\leq \tfrac{1}{\eta_i} \Et{\norm{\x_i^{t+1}-\x_i^t}^2} + \eta_i\Et{\norm{\g_i^t -\nabla f_i(\x_i^\star)}^2}\\
&\leq \tfrac{1}{\eta_i} \Et{\norm{\x_i^{t+1}-\x_i^t}^2} + 2A_i\eta_iD_{f_i}(\x_i^t,\x^\star_i) + B_i\eta_i \phi_i^t + 2\eta_i\sigma_i^2\label{i3}
\end{align}
where we have also used \eqref{a11}. 
Substituting \eqref{i1}, \eqref{i2}, and \eqref{i3}, we obtain
\begin{align}
\Et{\Psi^{t+1} + \sum_{i\in\V}\frac{1}{\eta_i}\norm{\x_i^{t+1}-\x_i^\star}^2}& \leq (\Psi^t + \sum_{i\in\V}\frac{1}{\eta_i}\norm{\x_i^t-\x_i^\star}^2)\\ 
&- 2\sum_{i\in\V}(1-A_i\eta_i)D_{f_i}(\x_i^t,\x_i^\star)\nonumber\\ &- 2\sum_{i\in\V}\Et{D_{h_i}(\x_i^{t+1},\x_i^\star)} \nonumber\\&+ \sum_{i\in\V}B_i\eta_i\phi_i^t + 2\eta_i\sigma_i^2\nonumber
\end{align}
Since $\varrho_i > 0$, the rate can be obtained by defining the Lyapunov function 
\begin{align}
\Upsilon^t = \Ex{\Psi^t + \sum_{i\in\V}\frac{1}{\eta_i}\norm{\x_i^t-\x_i^\star}^2 + \sum_{i\in\V}\frac{B_i\eta_i}{\varrho_i}\phi_i^t + 2\sum_{i\in\V}D_{h_i}(\x_i^t,\x_i^\star)}. \label{lyap}
\end{align}
which is always non-negative. From \eqref{w1}, we have that
\begin{align}
	\Upsilon^1 \leq 2EW_1(\rho+\frac{1}{\rho}) + W_1\sum_{i\in\V}\frac{1}{\eta_i} + \sum_{i\in\V}\frac{B_i\eta_i}{\varrho_i}\phi_i^1 + 2VW_1
\end{align}


Then adding \eqref{i1}, \eqref{i2}, and \eqref{i3} over all $i\in\V$, taking full expectation, and adding \eqref{a22} we obtain:
\begin{align}
\Upsilon^{t+1}-\Upsilon^t &\leq 2\sum_{i\in\V}(A_i\eta_i +\tfrac{B_iC_i\eta_i}{\varrho_i}-1)\Ex{D_{f_i}(\x_i^t,\x_i^\star)} \\
&-  2\sum_{i\in\V}\Ex{D_{h_i}(\x_i^{t},\x_i^\star)}+2\sum_{i\in\V}\sigma_i^2\eta_i \nonumber\\
&\leq -\sum_{i\in\V}\Ex{D_{f_i}(\x_i^t,\x_i^\star)} -  2\sum_{i\in\V}\Ex{D_{h_i}(\x_i^{t},\x_i^\star)}+2\sum_{i\in\V}\sigma_i^2\eta_i 
\end{align}
where we have assumed that $\frac{1}{2\eta_i} \geq A_i + \frac{B_iC_i}{\varrho_i}$ for all $i \in\V$. Re-arranging and using the fact that the Bregman divergences are non-negative, we obtain
\begin{align}
	\Delta^t = \sum_{i\in\V} \Ex{D_{f_i+h_i}(\x_i^t,\x_i^\star)} \leq \Upsilon^{t}-\Upsilon^{t+1} + 2\sum_{i\in\V}\sigma_i^2\eta_i
\end{align}
Finally, taking telescopic sum and using the fact that $\Upsilon^{T+1} \geq 0$, it can be seen that
\begin{align}
\frac{1}{T}\sum_{t=1}^T \Delta^t &\leq \frac{\Upsilon^1}{T} + 2\sum_{i\in\V}\sigma_i^2\eta_i \\
&\leq \left(2E\left(\rho+\frac{1}{\rho}\right) + \sum_{i\in\V} \frac{1}{\eta_i} + 2V\right)\frac{W_1}{T} + \sum_{i\in\V}\frac{B_i\eta_i\phi_i^1}{\varrho_i T} + 2\sum_{i\in\V} \sigma_i^2\eta_i 
\end{align}
which is the required result. 
\end{proof}

In the next section, we re-purpose existing centralized stochastic gradient methods to make them work in the decentralized setting. Lemma \ref{lema2} will be instrumental in studying these algorithms and each will correspond to a unique value of $A_i$, $B_i$, $C_i$, and $\varrho_i$. Observe here that the term $G_\rho$ appears separately and depends only on graph related quantities $V$ and $E$, as well as the parameter $\rho$. Generally, if we select a constant $\rho$, then for any graph, $G_\rho = \O(E)$.

\section{Algorithm Variants}\label{variants}
We present several variants of the proposed decentralized stochastic PD algorithm. Of particular interest are settings where $\{K_i\}$ are large and the full gradients $\nabla f_i(\x_i^t)$ cannot be evaluated at every iteration. Instead, each node has access to a first-order oracle that, for a given $\x_i^t$ and index $k_i^t$, outputs $\nabla f_i(\x_i^t,\xi_i^{k_i^t})$. Although complete expressions for all bounds are provided within the proofs, the main results will be stated using $O(\cdot)$ notation, which subsumes universal constants and the initialization dependent constant $W_1$. Further, we will assume that a constant $\rho$ is selected so that $G_\rho = \O(E)$. 

\subsection{Decentralized Stochastic Primal-Dual}
The simplest case, analogous to the classical SGD algorithm, is when $\g_i^t = \nabla f_i(\x_i^t,\xi_i^{k_i^t})$ where $k_i^t$ is selected uniformly at random from $\{1, \ldots, {K_i}\}$. With such a choice, it can be seen that $\g_i^t$ is an unbiased estimator of $\nabla f_i(\x_i^t)$, since 
\begin{align}
\Et{\g_i^t} = \mathbb{E}_{k_i^t}[\nabla f_i(\x_i^t,\xi_i^{k_i^t})] &= \frac{1}{{K_i}}\sum_{k=1}^{K_i} \nabla f_i(\x_i^t,\xi_i^k)= \nabla f_i(\x_i^t). \label{sgdunbiased}
\end{align}
where $\mathbb{E}_{k_i^t}[\cdot]$ denotes the expectation with respect to the random variable $k_i^t$. The stochastic PD is noise-dominated, and its performance depends on the gradient variance $\sigma^2 := (\sum_{i\in\V}\sigma_i)^2$, where
\begin{align}
\sigma_i^2 := \frac{1}{{K_i}}\sum_{k=1}^{K_i} \norm{\nabla f_i(\x_i^\star,\xi_i^{k}) - \nabla f_i(\x_i^\star)}^2
\end{align}
The following lemma provides a bound on the oracle complexity of online decentralized stochastic PD algorithm in terms of the cumulative Bregman divergence. 

\begin{lemma}\label{spd}
	Under Assumption \ref{a1}, the Decentralized Stochastic PD algorithm has an oracle complexity of $\O(\frac{\sigma^2}{\epsilon^2}+\frac{E}{\epsilon})$. The communication complexity is also the same. 
\end{lemma}
{ \begin{proof}[Proof of Lemma \ref{spd}]
We begin with verifying that Assumption \ref{a2} holds for the decentralized stochastic PD algorithm. We have already verified that $\g_i^t$ is unbiased [cf. \eqref{sgdunbiased}]. To see \eqref{a11}-\eqref{a22}, we have from the definition of $k_i^t$ that
\begin{align}
&\Et{\norm{\g_i^t - \nabla f_i(\x_i^\star)}^2}\\ & = \Et{\norm{\nabla f_i(\x_i^t,\xi_i^{k_i^t}) - \nabla f_i(\x_i^\star,\xi_i^{k_i^t}) + \nabla f_i(\x_i^\star,\xi_i^{k_i^t}) - \nabla f_i(\x_i^\star)}^2} \nonumber\\
&\hspace{-0.5cm}\leq 2\Et{\norm{\nabla f_i(\x_i^t,\xi_i^{k_i^t}) - \nabla f_i(\x_i^\star,\xi_i^{k_i^t})}^2} + 2\Et{\norm{\nabla f_i(\x_i^\star,\xi_i^{k_i^t}) - \nabla f_i(\x_i^\star)}^2} \label{osdpdpp} \\
& \hspace{-0.5cm}= \frac{2}{K_i}\sum_{k=1}^{K_i} \norm{\nabla f_i(\x_i^t,\xi_i^{k}) - \nabla f_i(\x_i^\star,\xi_i^{k})}^2 + \frac{2}{{K_i}}\sum_{k=1}^{K_i} \norm{\nabla f_i(\x_i^\star,\xi_i^{k}) - \nabla f_i(\x_i^\star)}^2 \label{osdpd2}
\end{align}
where \eqref{osdpdpp} follows from \eqref{expineq}. Bounding the first term in \eqref{osdpd2} from \eqref{dfbound} and recognizing that the second term is the variance of $\nabla f_i(\x_i^\star,\xi_i^{k_i^t})$, we obtain the constants required in Assumption \ref{a2} as $A_i = 2L_i$, $B_i = C_i = 0$, $\varrho_i = 1$, and 
\begin{align}
\sigma_i^2 =  \frac{1}{{K_i}}\sum_{k=1}^{K_i} \norm{\nabla f_i(\x_i^\star,\xi_i^{k}) - \nabla f_i(\x_i^\star)}^2. \label{sigma2}
\end{align} 
Therefore, we need to choose $\eta_i \leq \frac{1}{4L_i}$ and the bound in \eqref{mainbound} is written as
\begin{align}
	\frac{W_1G_\rho}{T} + \sum_{i\in \V}\left[\frac{W_1}{\eta_i T} + \sigma_i^2\eta_i\right] 
\end{align}
Setting $\eta_i = \frac{\sqrt{W_1}}{\sigma_i\sqrt{T}}$, we obtain a bound of $\O(\frac{E}{T} + \frac{\sigma}{\sqrt{T}})$ on the cumulative Bregman divergence where $\sigma := \sum_{i\in\V}\sigma_i$. Conversely, the number of oracle calls required to reach an $\epsilon$-optimal solution is given by $\O(\frac{\sigma^2}{\epsilon^2} + \frac{E}{\epsilon})$. Since a single round of communication, involving an exchange of primal and an exchange of dual variables, occurs at every $t$, the communication complexity is also the same. 
\end{proof}}

As in the classical SGD setting, the bound does not depend on the number of data points $K_i$, and can therefore be applied on the online setting where data arrives sequentially and indefinitely. In the online setting, the performance can be improved via mini-batching, where each node makes multiple oracle calls for the same $\x_i^t$ at each iteration.

\subsection{Decentralized SAGA Primal-Dual}
We first analyze the decentralized variant of the classical SAGA algorithm \cite{defazio2014saga}. In SAGA, variance reduction is achieved by maintaining $K_i$ gradient vectors corresponding to each data point. At each iteration, only one of these vectors is updated. The decentralized SAGA algorithm is summarized in Algorithm \ref{sagapd}. 

\begin{algorithm}
	\caption{Decentralized SAGA-PD}\label{sagapd}
	\begin{algorithmic}
		\STATE Initialize $\{\x_i^1\}$,  $\{\lam^1_{ij}\}$
		\STATE Set $\xt_i^1(k) = \x_i^1$ for all $1\leq k \leq K_i$ and $i\in\V$
		\FOR { $\textit{t} = 1,2,... $ } 
		\STATE pick  $k_i^t \hspace{2mm}\in\hspace{2mm} \{1,2,...,{K_i}\}$ uniformly and set	$\xt_i^{t+1}(k_i^t) = \x_i^t$ and $\xt_i^{t+1}(j) = \xt_i^t(j)$ for $j\neq k_i^t$.
		\STATE set $\g_i^t = \nabla f_i(\x_i^t,\xi_i^{k_i^t}) - \nabla f_i(\xt_i^t(k_i^t),\xi_i^{k_i^t})+ \frac{1}{K_i}\sum_{k=1}^{K_i} \nabla f_i(\xt_i^t(k),\xi_i^k)$
		\STATE update $\x_i^{t+1}$ using \eqref{primup} for all $i\in\V$,
		\STATE update $\lam_{ij}^{t+1}$ using \eqref{dualup} for all $i\in\V$, $j\in\N_i$
		\ENDFOR
	\end{algorithmic}
\end{algorithm}
Since $k_i^t$ is drawn uniformly at random, the gradient is unbiased. To see this, consider that
\begin{align}
\Et{\g_i^t} &= \frac{1}{K_i}\sum_{k=1}^{K_i} \nabla f_i(\x_i^t,\xi_i^k) - \frac{1}{K_i}\sum_{k=1}^{K_i}\nabla f_i(\xt_i^t(k),\xi_i^k) + \frac{1}{K_i}\sum_{k=1}^{K_i}\nabla f_i(\xt_i^t(k),\xi_i^{k})\nonumber
\end{align}
which equals $\nabla f_i(\x_i^t)$ since the last two terms cancel. The following lemma characterizes the iteration complexity of SAGA-PD. 
It can be seen that the oracle complexity of SAGA in Lemma \ref{saga-pd} matches that of the centralized version in \cite{defazio2014saga} for a single-node network.

\begin{lemma}\label{saga-pd}
	The SAGA PD algorithm incurs an  oracle complexity and communication complexity of  $\O\left(\frac{E+V\min\{L+K,\sqrt{LK}\}}{\epsilon}\right)$ for the choice $\eta = \frac{1}{2\sqrt{L\max\{K,16L\}}}$.
\end{lemma}
{ \begin{proof}[Proof of Lemma \ref{saga-pd}]
	 The proof for \eqref{a11}-\eqref{a22} proceeds by taking expectation with respect to $k_i^t$ and using \eqref{expineq} as follows:
	\begin{align}
	&\Et{\norm{\g_i^t - \nabla f_i(\x_i^\star)}^2}\nonumber \\
	&= \Et{\norm{\nabla f_i(\x_i^t,\xi_i^{k_i^t}) - \nabla f_i(\xt_i^t(k_i^t),\xi_i^{k_i^t}) + \frac{1}{K_i}\sum_{k=1}^{K_i} \nabla f_i(\xt_i^t(k),\xi_i^{k}) - \nabla f_i(\x_i^\star)}^2} \nonumber\\
	&\leq 2\Et{\norm{\nabla f_i(\x_i^t,\xi_i^{k_i^t}) -\nabla f_i(\x_i^\star,\xi_i^{k_i^t})}^2} + 2\Et{\norm{\nabla f_i(\xt_i^t,\xi_i^{k_i^t}) -\nabla f_i(\x_i^\star,\xi_i^{k_i^t})}^2} \nonumber
	\end{align}
	Here, the first term can be bounded using \eqref{dfbound} as
	\begin{align}
	\Et{\norm{\nabla f_i(\x_i^t,\xi_i^{k_i^t}) -\nabla f_i(\x_i^\star,\xi_i^{k_i^t})}^2} &= \frac{1}{K_i}\sum_{k=1}^{K_i}\norm{\nabla f_i(\x_i^t,\xi_i^{k}) -\nabla f_i(\x_i^\star,\xi_i^{k})}^2 \\
	&\leq 2L_iD_{f_i}(\x_i^t,\x_i^\star)
	\end{align}
	The second term can be written as
	\begin{align}
	\Et{\norm{\nabla f_i(\xt_i^t,\xi_i^{k_i^t}) -\nabla f_i(\x_i^\star,\xi_i^{k_i^t})}^2} &= \frac{1}{K_i}\sum_{k=1}^{K_i} \norm{\nabla f_i(\xt_i^t,\xi_i^{k}) -\nabla f_i(\x_i^\star,\xi_i^{k})}^2 =:\phi_i^t \nonumber
	\end{align}
	so that $A_i = 2L_i$, $B_i = 2$, and $\sigma_i^2 = 0$. To establish \eqref{a22}, observe that
	\begin{align}
	\phi_i^{t+1} = \frac{1}{K_i}\sum_{k=1}^{K_i} \norm{\nabla f_i(\xt_i^{t+1},\xi_i^{k}) -\nabla f_i(\x_i^\star,\xi_i^{k})}^2
	\end{align}
	where $\xt_i^{t+1}$ is $\x_i^t$ with probability $\frac{1}{K_i}$ and $\xt_i^t$ with probability $1-\frac{1}{K_i}$. Therefore, we obtain
	\begin{align}
	\phi_i^{t+1} &= \Et{ \frac{1}{K_i}\sum_{k=1}^{K_i} \norm{\nabla f_i(\xt_i^{t+1},\xi_i^{k}) -\nabla f_i(\x_i^\star,\xi_i^{k})}^2}\\
	&= \frac{1}{K_i}\sum_{k=1}^{K_i}\left[\frac{K_i-1}{K_i}\norm{\nabla f_i(\xt_i^t,\xi_i^{k}) -\nabla f_i(\x_i^\star,\xi_i^{k})}^2\right.\nonumber\\ &\left. + \frac{1}{K_i}\norm{\nabla f_i(\x_i^t,\xi_i^{k}) -\nabla f_i(\x_i^\star,\xi_i^{k})}^2\right]\nonumber\\
	&\leq \left(1-\frac{1}{K_i}\right)\phi_i^t + \frac{2L_i}{K_i}D_{f_i}(\x_i^t,\x_i^\star)
	\end{align}
	so that $\varrho_i = \frac{1}{K_i}$ and $C_i = \frac{L_i}{K_i}$. In this case, we also have that $\phi_i^1 \leq 2L_iD_{f_i}(\xt_i^1,\x_i^\star) = 2L_iD_{f_i}(\x_i^1,\x^\star)\leq 2L_iW_1$. Hence, Lemma \ref{lema2} requires us to choose $\eta_i \leq \frac{1}{8L_i}$, and the bound in \eqref{mainbound} becomes
	\begin{align}
		\frac{1}{T}\sum_{t=1}^T \Delta^t &\leq \frac{W_1G_\rho}{T} + \frac{W_1}{T}\sum_{i\in \V}\left[\frac{1}{\eta_i} + 4\eta_iL_iK_i\right] \\
		&= \frac{W_1}{T}\left(G_\rho + \frac{1}{2}\sum_{i\in\V} \min\{16L_i+K_i,8\sqrt{L_iK_i}\}\right)
	\end{align}
where we have used the optimal choice of $\eta_i$, given by $\eta_i = \frac{1}{2\sqrt{L_i\max\{K_i,16L_i\}}}$. Let us simplify the bound by setting $K_i = K$ and $L_i = L$ for all $i\in \V$, which yields
	\begin{align}
		\min_{1\leq t \leq T}\Delta_t \leq \O\left(\frac{E+V\min\{L+K,\sqrt{LK}\}}{T}\right)
	\end{align}
	and translates to an oracle complexity of $ \O\left(\frac{E+V\min\{L+K,\sqrt{LK}\}}{\epsilon}\right)$. Since a single round of communication, involving an exchange of primal and an exchange of dual variables between the neighbors, occurs at every $t$, the communication complexity is also the same. 
\end{proof}}

\subsection{Decentralized Stochastic Variance Reduced Primal-Dual}
Next, we present a variance reduction approach inspired from the SVRG++ algorithm \cite{allen2016improved} and summarized in Algorithm \ref{svrpd}. The algorithm proceeds in epochs indexed by $s$ and a complete gradient must be evaluated at every epoch. Additionally, two oracle calls are made at every iteration in order to calculate the gradient approximation; see Algorithm \ref{svrpd}. As in SVRG, the gradient approximation in SVR-PD is unbiased and its variance reduces with the number of epochs $S$. The following lemma provides the oracle complexity bound.

\begin{algorithm}[ht]
	\caption{Decentralized SVR-PD}\label{svrpd}
	\begin{algorithmic}
	\STATE Initialize $\{\xt^1_i = \x_i^0\}$,  $\{\lam^0_{ij}\}$ 
	\FOR { $\textit{s} = 1,2,..., S$ }  
		\FOR { $t=0,1,2,....., m_s-1$ }
				\STATE pick  $k_i^t \hspace{2mm}\in\hspace{2mm} \{1,2,...,{K_i}\}$ uniformly and set $\g_i^t = \nabla f_i(\x_i^t,\xi_i^{k_i^t}) -  \nabla f_i(\xt_i^s,\xi_i^{k_i^t}) + \nabla f_i(\xt_i^s)$
				\STATE update $\x_i^{t+1}$ using \eqref{primup} for all $i\in\V$,
				\STATE update $\lam_{ij}^{t+1}$ using \eqref{dualup}	for all $i\in\V$, $j\in\N_i$ 		
			\ENDFOR
		\STATE set $\xt_i^{s+1} = \frac{1}{m_s}\sum_{t=1}^{m_s} \x_i^t$ for $i\in\V$ 
	    \STATE for the next epoch, set $\x_i^0 = \x_i^m$ and $\lam_{ij}^0 = \lam_{ij}^m$ for all $i\in\V$ and $j\in\N_i$
	    \ENDFOR
	\end{algorithmic}
\end{algorithm}

\begin{lemma}\label{svr-pd}
	For constant $\eta = \frac{1}{6L}$ and $m_{s+1} = 2^sm_1$, the SVR-PD algorithm incurs an  oracle complexity of $\O(\frac{E+VL}{\epsilon} + K\log(\frac{E+VL}{\epsilon}))$ .  The communication complexity of SVR-PD is $\O(\frac{E+VL}{\epsilon})$. 
\end{lemma}

{ \begin{proof}[Proof of Lemma \ref{svr-pd}]
	Since $k_i^t$ is selected uniformly at random, it is clear that 
	\begin{align}
	\Et{\g_t^i} &= \frac{1}{{K_i}}\sum_{k=1}^{K_i}\left[\nabla f_i(\x_i^t,\xi_i^{k}) - \nabla f_i (\xt_i^s,\xi_i^{k})  \right] + \nabla f_i(\xt_i^s) = \nabla f_i(\x_i^t)
	\end{align}
	Next, to establish \eqref{a11}-\eqref{a22}, we introduce $\nabla f_i(\x_i^\star,\xi_i^{k_i^t})$ and expand as follows:
	\begin{small}
	\begin{gather}
	\Et{\norm{\g_i^t - \nabla f_i(\x_i^\star)}^2} = \Et{\norm{\nabla f_i(\x^t_i,\xi_i^{k_i^t}) - \nabla f_i(\xt_i^s,\xi_i^{k_i^t}) + \nabla f_i(\xt_i^s) - \nabla f_i(\x_i^\star)}^2} \\
	 =  \Et{\norm{\nabla f_i(\x_i^t,\xi_i^{k_i^t}) - \nabla f_i(\x_i^\star,\xi_i^{k_i^t}) +  \nabla f_i(\x_i^\star,\xi_i^{k_i^t}) - \nabla f_i(\xt_i^s,\xi_i^{k_i^t}) + \nabla f_i(\xt_i^s) - \nabla f_i(\x_i^\star)}^2} \nonumber\\
	\leq 2\Et{\norm{\nabla f_i(\x_i^t,\xi_i^{k_i^t}) - \nabla f_i(\x_i^\star,\xi_i^{k_i^t})}^2} + 2\Et{\norm{\nabla f_i(\x_i^\star,\xi_i^{k_i^t}) - \nabla f_i(\xt_i^s,\xi_i^{k_i^t})}^2} \label{svrpp}
	\end{gather}
	\end{small}
	where \eqref{svrpp} follows from \eqref{expineq}. The first term can be bounded using \eqref{dfbound} as
	\begin{align}
	\Et{\norm{\nabla f_i(\x_i^t,\xi_i^{k_i^t}) - \nabla f_i(\x_i^\star,\xi_i^{k_i^t})}^2} &= \frac{1}{{K_i}}\sum_{k=1}^{K_i}\norm{\nabla f_i(\x_i^t,\xi_i^k) - \nabla f_i(\x_i^\star,\xi_i^k)}^2\nonumber \\ &\leq 2L_iD_{f_i}(\x_i^t,\x_i^\star).
	\end{align}	
	For the second term of \eqref{svrpp}, we again use \eqref{dfbound} to obtain 
	\begin{align}
	\Et{\norm{\nabla f_i(\x_i^\star,\xi_i^{k_i^t}) - \nabla f_i(\xt_i^s,\xi_i^{k_i^t})}^2} & = \frac{1}{{K_i}}\sum_{k=1}^{K_i} \norm{\nabla f_i(\x_i^\star,\xi_i^k) - \nabla f_i(\xt_i^s,\xi_i^k)}^2\\
	& \leq 2L_iD_{f_i}(\xt_i^s,\x_i^\star) =:\phi_i^s\nonumber \label{svrphi}
	\end{align}
	Note that $\phi_i^s$ does not depend on $t$, implying that $\varrho_i = 0$ and $C_i = 0$. Substituting, we also obtain $A_i = 2L_i$,  $B_i = 2$, and $\sigma^2_i = 0$. 
	
	Since Lemma \ref{lema2} requires that $\varrho_i > 0$, the proof of Lemma \ref{lema2} cannot be directly used. Instead, let us define the function for iteration $t$ and epoch $s$:
	\begin{align}
	\Upsilon^t_s = \Ex{\Psi^t + \sum_{i\in\V}\frac{1}{\eta_i}\norm{\x_i^t-\x_i^\star}^2 + 2\sum_{i\in\V}D_{h_i}(\x_i^t,\x_i^\star)}
	\end{align}
	which yields
	\begin{align}
	\Upsilon^{t+1}_s - \Upsilon^t_s &\leq -2\sum_{i\in\V}(1-2L_i\eta_i)\Ex{D_{f_i}(\x_i^t,\x_i^\star)} -2\sum_{i\in\V}D_{h_i}(\x_i^t,\x_i^\star)+ 2\sum_{i\in\V}\eta_i\phi_i^s \\
	&\leq -2\sum_{i\in\V}(1-2L_i\eta_i)\Ex{D_{f_i+h_i}(\x_i^t,\x_i^\star)} + 4\sum_{i\in\V}\eta_iL_i D_{f_i+h_i}(\xt_i^s,\x_i^\star).
	\end{align}
	In order to solve this recursion, let us simplify the notation by denoting $\Delta_s^t := \sum_{i}D_{f_i+h_i}(\x_i^t,\x_i^\star)$ and $\Dt_s:= \sum_{i}D_{f_i+h_i}(\xt_i^s,\x_i^\star)$. Also select $\eta_i = \frac{1}{6L_i}$ so that $\eta_iL_i = \frac{1}{6}$, which yields
	\begin{align}
	\Upsilon^{t+1}_s - \Upsilon^t_s &\leq -2(1-\tfrac{1}{3})\Delta_s^t + \tfrac{2}{3}\Dt_s \\
	\Rightarrow \tfrac{4}{3}\Delta_s^t &\leq \tfrac{2}{3}\Dt_s + \Upsilon^t_s - \Upsilon^{t+1}_s.   
	\end{align}
	Rearranging, summing over $t = 0, \ldots, m_s-1$, and using the convexity of $D_{f_i+h_i}(\x,\x_i^\star)$ in its first argument, we obtain
	\begin{align}
	\tfrac{4}{3}\Dt_{s+1} \leq \frac{4}{3m_s}\sum_{t=0}^{m_s-1}\Delta_s^t \leq \tfrac{2}{3}\Dt_s + \frac{\Upsilon_s^0-\Upsilon_s^{m_s}}{m_s} \label{svrgtel}
	\end{align}
	Since $\Upsilon_s^{m_s} = \Upsilon^0_{s+1}$, we can rewrite \eqref{svrgtel}  as
	\begin{align}
	\Dt_{s+1} + \frac{3\Upsilon_{s+1}^0}{2m_{s+1}} \leq \frac{1}{2}\left(\Dt_s + \frac{3\Upsilon_s^0}{2m_{s}}\right)
	\end{align}
	where we have used the fact that $m_{s+1} = 2m_s$. Taking telescopic sum, we obtain 
	\begin{align}
	\Dt_{S+1} + \frac{3\Upsilon_{S+1}^0}{2m_{S+1}} \leq 2^{-S}(\Dt_1 + \frac{3\Upsilon_1^0}{2m_1}).
	\end{align}
	Here, $\Dt_1 \leq 2VW_1$ and $\Upsilon_1^0 \leq W_1(G_\rho + 6\sum_i L_i)$, so that
	\begin{align}
		\Dt_{S+1} \leq 2^{-S}W_1\left(2V+\frac{3(G_\rho+6\sum_i L_i)}{2m_1}\right)
	\end{align}
	To simplify the expression, let us put $L_i = L$ for all $i\in\V$, which yields
	\begin{align}
	\Dt_{S+1} \leq \O(2^{-S}(E+VL))
	\end{align}
	where the $\O(\cdot)$ notation subsumes initialization-dependent constants like $W_1$ and $m_1$.
	
	Therefore in order to obtain $\O(\epsilon)$ cumulative Bregman divergence, we require $\O(\log(\frac{E+VL}{\epsilon}))$ epochs. Further, the total number of oracle calls in $S$ epochs is $\O(2^S + SK)$. Therefore, the required oracle complexity becomes $\O(\frac{E+VL}{\epsilon} + K\log(\frac{E+VL}{\epsilon}))$. Also, there are $m_s$ rounds of communication at epoch $s$, which translates to a total of $\sum_s m_s = \O(2^S)$ rounds of communications, each involving an exchange of primal and an exchange of dual variables. The communication complexity is therefore $\O(\frac{E+VL}{\epsilon})$.
\end{proof}}

The iteration complexity achieved by the SVR-PD is the same as that of the SVRG++ algorithm \cite{allen2016improved}. It can be seen that for $K \ll \frac{1}{\epsilon}$, SVR-PD is $\O(\frac{1}{\epsilon})$ times faster than decentralized stochastic PD algorithm. Further, the SVR-PD algorithm does not depend on the data variance $\sigma^2$, but only on the smoothness parameter $L$. In general, SVR-PD enjoys all the benefits of SVRG++ in the decentralized setting.

\subsection{Decentralized Loopless SVR-PD}
We present the loopless version of SVR-PD along the lines of L-SVRG in \cite{kovalev2019don}. In the literature, the L-SVRG algorithm has only been analyzed for the strongly convex case. Therefore, the present analysis also provides a new bound for L-SVRG for general convex objectives. 

\begin{algorithm}[ht]
	\caption{Decentralized LSVR-PD}\label{lsvrpd}
	\begin{algorithmic}
		\STATE Initialize $\{\xt^1_i = \x_i^0\}$,  $\{\lam^0_{ij}\}$
		\FOR { $\textit{t} = 1,2,... $ }  
		\STATE pick  $k_i^t \hspace{2mm}\in\hspace{2mm} \{1,2,...,{K_i}\}$ uniformly and set
		\begin{align}
		\g_i^t = \nabla f_i(\x_i^t,\xi_i^{k_i^t}) -  \nabla f_i(\xt_i^t,\xi_i^{k_i^t}) + \nabla f_i(\xt_i^t)
		\end{align}
		\STATE update $\x_i^{t+1}$ using \eqref{primup} for all $i\in\V$,
		\STATE update $\lam_{ij}^{t+1}$ using \eqref{dualup} for all $i\in\V$, $j\in\N_i$
		\STATE pick $\Xi_i^t \sim $Bernoilli$(\frac{1}{K_i})$ and set $\xt_i^{t+1} = \Xi_i^t\x_i^t + (1-\Xi_i^t)\xt_i^t$
		\ENDFOR
	\end{algorithmic}
\end{algorithm}

The following theorem characterizes the oracle complexity. Observe that the complexity is almost the same as that of SAGA-PD but worse than SVR-PD. 

\begin{lemma}\label{lsvr-pd}
	The LSVR-PD algorithm incurs an  oracle complexity of incurs an  oracle complexity and communication complexity of  $\O\left(\frac{E+V\min\{L+K,\sqrt{LK}\}}{\epsilon}\right)$ for the choice $\eta = \frac{1}{\sqrt{2L\max\{K,32L\}}}$.
\end{lemma}
\begin{proof}[Proof of Lemma \ref{lsvr-pd}]
	It can be seen that $\Et{\g_i^t} = \nabla f_i(\x_i^t)$. In order to obtain \eqref{a22}, we take expectation with respect to the independent variables $k_i^t$ and $\Xi_i^t$ and introduce $\nabla f_i(\x_i^\star,\xi_i^{k_i^t})$ as follows:
	\begin{small}
	\begin{align}
	&\Et{\norm{\g_i^t-\nabla f_i(\x_i^\star)}^2} =  \Et{\norm{\nabla f_i(\x^t_i,\xi_i^{k_i^t}) - \nabla f_i(\xt_i^t,\xi_i^{k_i^t}) + \nabla f_i(\xt_i^t) - \nabla f_i(\x_i^\star)}^2} \nonumber\\
	& =  \Et{\norm{\nabla f_i(\x_i^t,\xi_i^{k_i^t}) - \nabla f_i(\x_i^\star,\xi_i^{k_i^t}) +  \nabla f_i(\x_i^\star,\xi_i^{k_i^t}) - \nabla f_i(\xt_i^t,\xi_i^{k_i^t}) + \nabla f_i(\xt_i^t) - \nabla f_i(\x_i^\star)}^2} \nonumber\\
	&\leq 2\Et{\norm{\nabla f_i(\x_i^t,\xi_i^{k_i^t}) - \nabla f_i(\x_i^\star,\xi_i^{k_i^t})}^2} + 2\Et{\norm{\nabla f_i(\x_i^\star,\xi_i^{k_i^t}) - \nabla f_i(\xt_i^t,\xi_i^{k_i^t})}^2} \label{lsvrpp}
	\end{align}
	\end{small}
	where \eqref{lsvrpp} follows from \eqref{expineq}. Proceeding as in the proof of Lemma \ref{svr-pd}, we obtain the bound:
	\begin{align}
	\Et{\norm{\g_i^t-\nabla f_i(\x_i^\star)}^2} \leq 4L_iD_{f_i}(\x_i^t,\x_i^\star) + 4L_iD_{f_i}(\xt_i^t,\x_i^\star).
	\end{align}
	Comparing with \eqref{a11}, we have that $A_i = 2L_i$, $B_i = 2L_i$, and $\phi_i^t := D_{f_i}(\xt_i^t,\x_i^\star)$. Further, from the update equation, it is clear that 
	\begin{align}
	\Ex{\phi_i^{t+1}} &= \Ex{D_{f_i}(\xt_i^{t+1},\x_i^\star)} = (1-\tfrac{1}{K_i})\Ex{D_{f_i}(\xt_i^t,\x_i^\star)} + \tfrac{1}{K_i}\Ex{D_{f_i}(\x_i^t,\x_i^\star)} \\
	&= (1-\tfrac{1}{K_i})\phi_i^t + \tfrac{1}{K_i}\Ex{D_{f_i}(\x_i^t,\x_i^\star)}
	\end{align}
	which implies that $\varrho_i = 1/K_i$ and $C_i = \frac{1}{K_i}$. From Lemma \ref{lema1}, it is therefore required to choose $\eta_i \leq \frac{1}{8L_i}$, and the resulting bound takes the form
	\begin{align}
		\min_{1\leq t \leq T}\Delta^t &\leq  \frac{W_1G_\rho}{T} + \frac{W_1}{T}\sum_{i\in \V}\left[\frac{1}{\eta_i} + 2L_iK_i\eta_i\right] \\
		& = \frac{W_1}{T}\left(G_\rho + \frac{1}{4}\sum_{i\in\V} \min\{32L_i+K_i,8\sqrt{2L_iK_i}\}\right)
	\end{align}
	where we have used the optimal choice of $\eta_i$, given by $\eta_i = \frac{1}{\sqrt{2L_i\max\{K_i,32L_i\}}}$. Simplifying this bound by setting $K_i = K$ and $L_i = L$ for all $i\in \V$, we obtain
	\begin{align}
		\min_{1\leq t \leq T}\Delta_t \leq \O\left(\frac{E+V\min\{L+K,\sqrt{LK}\}}{T}\right)
	\end{align}
	which translates to an oracle complexity of $\O\left(\frac{E+V\min\{L+K,\sqrt{LK}\}}{\epsilon}\right)$. Since one round of communication, involving an exchange of primal and an exchange of dual variables, occurs at every $t$, the communication complexity is also the same. 
\end{proof}
\section{Accelerated SVR-PD}\label{accsec}
Algorithm \ref{asvrpdmm} puts forth the Accelerated SVR-PD  (ASVR-PD) Algorithm. The proposed algorithm is partly inspired from the accelerated SVRG algorithm proposed in \cite{shang2018asvrg}. However, the adaptation is not straightforward, as the decentralized algorithm still requires more message passing. As in the SVR-PD algorithm, the proposed ASVR-PD algorithm also proceeds in epochs of variable sizes. In addition, more message passing is achieved by making an oracle call with some probability $p_s$ only. In the remaining iterations, the nodes simply use the stored gradient, but still exchange messages with each other. To ensure that the resulting gradient is unbiased, we select it as 
\begin{align}
\g^t_i = \begin{cases} \frac{1}{p_s}\left(\nabla f_i(\x^t_i,\xi_i^{k_i^t}) -  \nabla f_i(\xt^s_i,\xi_i^{k_i^t})\right) + \nabla f_i(\xt^s_i), \quad \text{with probability $p_s$} \\
\nabla f_i(\xt^s_i), \quad \text{with probability $1-p_s$}
\end{cases}
\end{align}
As earlier, $\Xi_i^t$ be the random variable that is 1 with probability $p_s$ and 0 with probability $1-p_s$. Then the gradient can be written as
\begin{align}
\g_i^t = \frac{\Xi_i^t}{p_s}\left(\nabla f_i(\x^t_i,\xi_i^{k_i^t}) -  \nabla f_i(\xt^s_i,\xi_i^{k_i^t})\right) + \nabla f_i(\xt^s_i)
\end{align}
Since $\Et{\Xi_{t}} = p_s$ and $\xi_i^{k_i^t}$ is independent of $\Xi_{t}$, we have that $\Et{\g_i^t} = \nabla f(\x_i^t)$. 

Additionally, the negative momentum parameter $\beta_s$ must be selected in a specific manner in order to utilize the telescopic sum; such a choice is typical of accelerated algorithms. The full algorithm is summarized in Algorithm \ref{asvrpdmm}.

\begin{algorithm}[ht]
	\caption{Accelerated SVR-PD}\label{asvrpdmm}
	\begin{algorithmic}
	\STATE Initialize $\xt^1_i = \yt^1_i$ for $i\in\V$, $\beta_1 = 1/2$ 
	\FOR { $\textit{s} = 1,2,... $ } 
		\STATE Set $\x^0_i = (1-\beta_s)\xt^s_i + \beta_s \yt^s_i$ and set $\y^0_i = \yt^s_i$
		\FOR { $t=0,1,2,....., m_s-1$ }			
			\STATE pick  $k_i^t \hspace{2mm}\in\hspace{2mm} \{1,2,...,{K_i}\}$ uniformly at random and set for each $i\in\V$
			\begin{align*}
			\g^t_i &= \begin{cases} \frac{1}{p_s}\left(\nabla f_i(\x^t_i,\xi_i^{k_i^t}) -  \nabla f_i(\xt^s_i,\xi_i^{k_i^t})\right) + \nabla f_i(\xt^s_i), & \text{with probability $p_s$} \\
			 \nabla f_i(\xt^s_i), & \text{with probability $1-p_s$}
			\end{cases}
			\end{align*}
			\STATE Update $\y_i^{t+1} =\pr_{\gamma_i h_i}\Big(\frac{\gamma_i}{\eta_i}[\y_i^t-\eta_i\g_i^t] + \gamma_i\sum_{j\in\N_i}[\rho\y_j^t-\lam_{ji}^t] \Big)$	for $i\in\V$.
			\STATE Update $\lam_{ij}^{t+1} = - \lam_{ji}^t + \rho (\y_j^t - \y_i^{t+1})$  for $i\in\V$
			\STATE Update $\x_i^{t+1} = (1-\beta_s)\xt_i^s + \beta_s \y_i^{t+1}$ for $i\in\V$			
		\ENDFOR
		\STATE Set $\xt_i^{s+1} = \frac{1}{m_s}\sum_{t=1}^{m_s} \x_i^t$  and $\yt_i^{s+1} = \y_i^{m_s}$ for $i\in\V$
		\STATE Set $\beta_{s+1} = \frac{\sqrt{\beta_s^4 + 4\beta_s^2} - \beta_s^2}{2}$ 
	\ENDFOR
	\end{algorithmic}
\end{algorithm}

Here, the $\y_i$-update can also be written as
\begin{align}
\y_i^{t+1} = \pr_{\eta_ih_i}\Big(\y_i^t - \eta \g_i^t + \eta\sum_{j\in\N_i}\lam_{ij}^{t+1}\Big)
\end{align}
The update equations imply that $\x_i^{t+1}-\x_i^t = \beta_s(\y_i^{t+1}-\y_i^t)$. Note that the number of samples used on average per epoch is $m_sp_s$ and the total number of samples used is $\sum_{s=1}^S m_sp_s$. 

\begin{lemma}\label{asvrlem}
	For $p_s = \beta_s$ and $m_s = \frac{K}{\beta_s}$, the oracle complexity of the Accelerated SVR-PD algorithm is $\O\left(\frac{\sqrt{K(E+V(K+L))}}{\sqrt{\epsilon}}\right)$. The communication complexity is $\O(\frac{E+V(K+L)}{\epsilon})$. 
\end{lemma}

{\begin{lemma}
	The gradient approximation $\g_i^t$ satisfies:
	\begin{align}
	\Et{\norm{\g_i^t - \nabla f_i(\x_i^t)}^2} \leq \frac{2L_i}{p_s}D_{f_i}(\xt_i^s,\x_i^t) \label{asvrpdmm-gt}
	\end{align}
	for all $0 \leq t \leq m_s-1$ and epoch $s$. 
\end{lemma}
\begin{proof}
	For the sake of brevity, let us denote:
	\begin{align}
	\sX_t  &= \nabla f_i(\x^t_i,\xi_i^{k_i^t}) & \Rightarrow \mx &= \nabla f_i(\x^t_i)\\
	\sY_t &= \nabla f_i(\xt^s_i,\xi_i^{k_i^t}) & \Rightarrow \my &= \nabla f_i(\xt^s_i)
	\end{align}
	and $\g_i^t = \frac{\Xi_t}{p_s}(\sX_t - \sY_t) + \my$ with $\Et{\g_i^t} = \mx$. Therefore, we have that
	\begin{align}
	\Et{\norm{\g_i^t-\mx}^2} &= p_s\Et{\norm{\frac{\sX_t - \sY_t}{p_s} + \my - \mx}^2} + (1-p_s)\norm{\my-\mx}^2\\
	&\hspace{-2cm}=p_s\left(\frac{1}{p_s^2}\Et{\norm{\sX_t-\sY_t}^2} + \norm{\mx-\my}^2 - \frac{2}{p_s}\norm{\mx-\my}^2\right) + (1-p_s)\norm{\mx-\my}^2 \nonumber\\
	& = \frac{1}{p_s}\Et{\norm{\sX_t-\sY_t}^2} - \norm{\mx-\my}^2
	\end{align}
	Dropping the last term and substituting, we obtain
	\begin{align}
	\Et{\norm{\g_i^t - \nabla f_i(\x_i^t)}^2} &\leq \frac{1}{p_s}\Et{\norm{\nabla f_i(\x_i^t,\xi_i^{k_i^t}) -  \nabla f_i(\xt_i^s,\xi_i^{k_i^t})}^2} \\
	&= \frac{1}{{K_i}p_s}\sum_{k=1}^{K_i} \norm{\nabla f_i(\x_i^t,\xi_i^k) -  \nabla f_i(\xt_i^s,\xi_i^k)}^2 \\
	& \leq \frac{1}{{K_i}p_s}\sum_{k=1}^{K_i} 2L_i D_{f_i}(\xt_i^s,\x_i^t) =  \frac{2L_i}{p_s} D_{f_i}(\xt_i^s,\x_i^t)
	\end{align}
	where we have used the inequality $\norm{\nabla f_i(\x)-\nabla f_i(\y)}^2 \leq 2 L_i D_{f_i}(\x,\y)$.
\end{proof}

We also establish the following preliminary equality.
\begin{lemma}
	It holds that
	\begin{align}
	\ip{\g_i^t}{\y_i^{t+1}-\y_i^t} + \tfrac{\beta_s}{2\eta_i}\norm{\y_i^{t+1}-\y_i^t}^2 &= \ip{\g_i^t}{\x_i^\star-\y_i^t} +  \tfrac{\beta_s}{2\eta_i}\norm{\y_i^t-\x_i^\star}^2 \nonumber \\- \tfrac{\beta_s}{2\eta_i}\norm{\y_i^{t+1}-\x_i^\star}^2 
	&+ \ip{\sum_{j\in\N_i}\lam_{ij}^{t+1}-\v_i^{t+1}}{\y_i^{t+1} - \x_i^\star} \label{prelim2}
	\end{align}
\end{lemma}
\begin{proof}
	Since  $\v_i^{t+1} + \g_i^t - \sum_{j\in\N_i}\lam_{ij}^{t+1} = \frac{1}{\eta_i}(\y_i^{t}-\y_i^{t+1})$ for some $\v_i^{t+1} \in \partial h_i(\y_i^{t+1})$, we have that
	\begin{align}
	2\eta_i\ip{\v_i^{t+1} + \g_i^t - \sum_{j\in\N_i}\lam_{ij}^{t+1}}{\y_i^{t+1}-\x_i^\star} &=  2\ip{\y_i^t-\y_i^{t+1}}{\y_i^{t+1}-\x_i^\star} \\
	& \hspace{-3.5cm}= \norm{\y_i^t-\x_i^\star}^2 - \norm{\y_i^{t+1}-\x_i^\star}^2 - \norm{\y_i^{t+1}-\y_i^t}^2 \nonumber\\
	\Rightarrow \ip{\g_i^t - \sum_{j\in\N_i}\lam_{ij}^{t+1}}{\y_i^{t+1}-\y_i^t} + \tfrac{1}{2\eta_i}\norm{\y_i^{t+1}-\y_i^t}^2 &= \ip{\g_i^t - \sum_{j\in\N_i}\lam_{ij}^{t+1}}{\x_i^\star-\y_i^t} \nonumber\\
	&\hspace{-3.5cm}+  \tfrac{1}{2\eta_i}\norm{\y_i^t-\x_i^\star}^2 - \tfrac{1}{2\eta_i}\norm{\y_i^{t+1}-\x_i^\star}^2 +  \ip{\v_i^{t+1}}{\x_i^\star - \y_i^{t+1}} \\
	\Rightarrow \ip{\g_i^t}{\y_i^{t+1}-\y_i^t} + \tfrac{1}{2\eta_i}\norm{\y_i^{t+1}-\y_i^t}^2 &\leq \ip{\g_i^t}{\x_i^\star-\y_i^t} +  \tfrac{1}{2\eta_i}\norm{\y_i^t-\x_i^\star}^2 \nonumber \\
	&\hspace{-3.5cm}- \tfrac{1}{2\eta_i}\norm{\y_i^{t+1}-\x_i^\star}^2 +\ip{\sum_{j\in\N_i}\lam_{ij}^{t+1}-\v_i^{t+1}}{\y_i^{t+1} - \x_i^\star}  
	\end{align}
\end{proof}

We are now ready to present the proof of Lemma \ref{asvrlem}. To this end, we redefine $\Psi^t$ as follows:
\begin{align}
\Psi_t^s :=  \frac{1}{2\rho}\sum_{i\in\V}\sum_{j\in\N_i}\norm{\rho(\y_i^t-\x_i^\star) - (\lam_{ij}^t-\lam_{ij}^\star)}^2 \label{bks}
\end{align} 
for any $t$ within an epoch $s$.

\begin{proof}
	We begin with the quadratic upper bound for $f_i$ for each $i \in \V$
	\begin{align}
	f_i(\x_i^{t+1}) &\leq f_i(\x_i^t) + \ip{\nabla f_i(\x_i^t)}{\x_i^{t+1}-\x_i^t} + \frac{L_i}{2}\norm{\x_i^t-\x_i^{t+1}}^2 \nonumber\\
	& = f_i(\x_i^t) + \ip{\g_i^t}{\x_i^{t+1}-\x_i^t} + \ip{\nabla f_i(\x_i^t)-\g_i^t}{\x_i^{t+1}-\x_i^t} + \frac{L_i}{2}\norm{\x_i^t-\x_i^{t+1}}^2 \\
	&\leq f_i(\x^t) + \ip{\g_i^t}{\x_i^{t+1}-\x_i^t} + \frac{1}{2\omega_i^s}\norm{\nabla f_i(\x_i^t)-\g_i^t}^2  + \frac{\omega_i^s+L_i}{2}\norm{\x_i^t-\x_i^{t+1}}^2 
	\end{align}
	where we have used the Peter-Paul inequality with parameter $\omega_i^s$. Recall that we have:
	\begin{align}
	\x_i^{t+1}-\x_i^t &= \beta_s(\y_i^{t+1}-\y_i^t) \label{moment2}
	\end{align}
	Substituting \eqref{moment2} and \eqref{asvrpdmm-gt}, and taking expectation, we obtain
	\begin{align}
	\Et{f_i(\x_i^{t+1}) - f_i(\x_i^t)} &\leq  \beta_s\Et{\ip{\g_i^t}{\y_i^{t+1}-\y_i^t}} + \frac{L_i}{\omega_i^s p_s}D_{f_i}(\xt_i^s,\x_i^t)  \\&+ \frac{\beta_s^2(\omega_i^s+L_i)}{2}\Et{\norm{\y_i^t-\y_i^{t+1}}^2}\nonumber\\
	&\hspace{-2cm}\leq \beta_s\Et{\ip{\g_i^t}{\y_i^{t+1}-\y_i^t}} + \frac{L_i}{\omega_i^s p_s}D_{f_i}(\xt_i^s,\x_i^t) + \frac{\beta_s}{2\eta_i}\Et{\norm{\y_i^t-\y_i^{t+1}}^2}
	\end{align}
	where the last inequality follows when $\eta_i \leq \frac{1}{\beta_s(\omega_i^s+L_i)}$. It follows from \eqref{prelim2} that
	\begin{align}
	\Et{f_i(\x_i^{t+1}) - f_i(\x_i^t)} &\leq  \beta_s\Et{\ip{\g_i^t}{\x_i^\star-\y_i^t}} + \frac{L_i}{\omega_i^sp_s}D_{f_i}(\xt_i^s,\x_i^t)  + \frac{\beta_s}{2\eta_i}\left(\norm{\x_i^\star-\y_i^{t}}^2 \right. \nonumber\\
	&\left.- \Et{\norm{\x_i^\star-\y_i^{t+1}}^2}\right)  + \beta_s\ip{\sum_{j\in\N_i}\lam_{ij}^{t+1}-\v_i^{t+1}}{\y_i^{t+1} - \x_i^\star}\\
	&\hspace{-2cm}=  \ip{\nabla f_i(\x_i^t)}{\beta_s\x_i^\star + (1-\beta_s)\xt_i^s - \x_i^t} + \frac{L_i}{\omega_i^sp_s}D_{f_i}(\xt_i^s,\x_i^t)  + \frac{\beta_s}{2\eta_i}\left(\norm{\x_i^\star-\y_i^{t}}^2\right. \nonumber\\
	&\left. - \Et{\norm{\x_i^\star-\y_i^{t+1}}^2}\right)  + \beta_s\ip{\sum_{j\in\N_i}\lam_{ij}^{t+1}-\v_i^{t+1}}{\y_i^{t+1} - \x_i^\star} 
	\end{align}
	Next, let us analyze the first two terms on the right. Expanding, and from the convexity of $f_i$, we have that 
	\begin{align}
	&\ip{\nabla f_i(\x_i^t)}{\beta_s\x_i^\star + (1-\beta_s)\xt_i^s - \x_i^t} + \frac{L_i}{\omega_i^sp_s}D_{f_i}(\xt_i^s,\x_i^t) \\
	& = \ip{\nabla f_i(\x_i^t)}{\beta_s\x_i^\star + (1-\beta_s - \tfrac{L_i}{\omega_i^sp_s})\xt_i^s + \tfrac{L_i}{\omega_i^sp_s}\x_i^t - \x_i^t} + \tfrac{L_i}{\omega_i^s}(f_i(\xt_i^s) - f_i(\x_i^t)) \nonumber\\
	& \leq f_i(\beta_s\x_i^\star + (1-\beta_s - \tfrac{L_i}{\omega_i^sp_s})\xt_i^s + \tfrac{L_i}{\omega_i^s}\x_i^t) - f_i(\x_i^t) + \tfrac{L_i}{\omega_i^sp_s}(f_i(\xt_i^s) - f_i(\x_i^t)) \\
	&\leq \beta_sf_i(\x_i^\star) + (1-\beta_s) f_i(\xt_i^s) - f_i(\x_i^t)
	\end{align}
	provided that $\beta_s + \frac{L_i}{\omega_i^sp_s} \leq 1$ or $\omega_i^s \geq \frac{L_i}{p_s(1-\beta_s)}$. 
	Combining, it can be seen that
	\begin{align}
	\Et{f_i(\x_i^{t+1}) - f_i(\x_i^\star)} &\leq (1-\beta_s)(f_i(\xt_i^s)-f_i(\x_i^\star)) +  \frac{\beta_s}{2\eta_i}\left(\norm{\x_i^\star-\y_i^{t}}^2\right. \nonumber\\
	&\left.- \Et{\norm{\x_i^\star-\y_i^{t+1}}^2}\right) +\beta_s\ip{\sum_{j\in\N_i}\lam_{ij}^{t+1}-\v_i^{t+1}}{\y_i^{t+1} - \x_i^\star} 
	\end{align}
	Summing over $i\in\V$ and using Lemma \ref{lema1}, we obtain
	\begin{align}
	\sum_{i\in\V}\Et{f_i(\x_i^{t+1}) - f_i(\x_i^\star)} &\leq (1-\beta_s)\sum_{i\in\V}(f_i(\xt_i^s)-f_i(\x_i^\star)) +  \sum_{i\in\V}\frac{\beta_s}{2\eta_i}\left(\norm{\x_i^\star-\y_i^{t}}^2- \right.\nonumber\\
	&\hspace{-3.2cm}\left. \Et{\norm{\x_i^\star-\y_i^{t+1}}^2}\right) + \tfrac{\beta_s}{2}(\Psi_t^s-\Et{\Psi_{t+1}^s}) + \beta_s\sum_{i\in\V}\ip{\sum_{j\in\N_i}\lam_{ij}^{\star}-\v_i^{t+1}}{\y_i^{t+1} - \x_i^\star} \label{precontr}
	\end{align}
	The last term can be bounded as follows. 
	\begin{align}
	\beta_s\ip{\sum_{j\in\N_i}\lam_{ij}^{\star}-\v_i^{t+1}}{\y_i^{t+1} - \x_i^\star} &= \beta_s\ip{\nabla f(\x_i^\star)+\v_i^\star}{\y_i^{t+1} - \x_i^\star} - \beta_s\ip{\v_i^{t+1}}{\y_i^{t+1} - \x_i^\star}\nonumber\\
	& \hspace{-1.5cm}\leq  \ip{\nabla f(\x_i^\star)+\v_i^\star}{\beta_s\y_i^{t+1} - \beta_s\x_i^\star} - \beta_s(h_i(\y_i^{t+1})-h_i(\x_i^\star))  \\
	& \hspace{-1.5cm} \leq  \ip{\nabla f(\x_i^\star)+\v_i^\star}{\x_i^{t+1} - \x_i^\star} - (1-\beta_s)\ip{\nabla f(\x_i^\star)+\v_i^\star}{\xt_i^s - \x_i^\star}\\
	&\hspace{-1.5cm}  - (h(\x_i^{t+1})  - h_i(\x_i^\star)) + (1-\beta_s)(h_i(\xt_i^s) - h_i(\x_i^\star))\nonumber
	\end{align}
	where we have used the convexity of $h_i$ which implies that $h_i(\x_i^{t+1}) \leq (1-\beta_s)h_i(\xt_i^s) + \beta_sh_i(\y_i^{t+1})$.
	
	Assuming that $\x^t$ collects all $\{\x_i^t\}$ and likewise $\xt^s$ collects all $\{\xt_i^s\}$, define the convex and non-negative function
	\begin{align}
	\Gamma(\x) &= \sum_{i\in\V}\Ex{D_{f_i+h_i}(\x_i,\x_i^\star)} 
	\end{align}
	Taking full expectation in \eqref{precontr}, it can be seen that
	\begin{align}
	\Gamma(\x^{t+1}) &\leq (1-\beta_s)\Gamma(\xt^s) +\sum_{i\in\V}\frac{\beta_s}{2\eta_i}(\Ex{\norm{\x_i^\star-\y^t_i}^2}-\Ex{\norm{\x_i^\star-\y^{t+1}_i}^2})\\ &+\frac{\beta_s}{2}\Ex{\Psi_t^s-\Psi_{t+1}^s}\nonumber
	\end{align}
	Summing over $t = 1, \ldots, m_s$, and from the convexity of $\Gamma$, we obtain
	\begin{align}
	\Gamma(\xt^{s+1}) \leq (1-\beta_s)\Gamma(\xt^s) +& \frac{\beta_s}{ 2m_s}\sum_{i\in\V}\frac{1}{\eta_i}\left(\Ex{\norm{\x_i^\star-\yt^s_i}^2}-\Ex{\norm{\x_i^\star-\yt^{s+1}_i}^2}\right)\\ &+\frac{\beta_s}{2m_s}\Ex{\Psi^s_1-\Psi^{s+1}_1}\nonumber
	\end{align}
	where we have also used the fact that $\yt^s = \y_i^0$ and $\yt^{s+1} = \y_i^m$. 	Suppose that we choose $m_s = m/\beta_s$ for some constant $m \geq 1$, then we have that
	\begin{align}
	\Gamma(\xt^{s+1}) \leq (1-\beta_s)\Gamma(\xt^s) + &\frac{\beta_s^2}{2 m}\sum_{i\in\V}\frac{1}{\eta_i}\left(\Ex{\norm{\x_i^\star-\yt^s_i}^2}-\Ex{\norm{\x_i^\star-\yt^{s+1}_i}^2}\right)\\ &+\frac{\beta_s^2}{2m}\Ex{\Psi^s_1-\Psi^{s+1}_1}\nonumber
	\end{align}
	If we choose $\frac{1}{\beta_{s-1}^2} = \frac{1-\beta_s}{\beta_s^2}$, then we obtain
	\begin{align}
	\sum_{i\in\V}\Ex{D_{f_i+h_i}(\xt_i^{S+1},\x_i^\star)} &\leq \frac{\beta_S^2}{2m} \left[\sum_{i\in\V}\Ex{mD_{f_i+h_i}(\xt_i^{1},\x_i^\star) + \frac{1}{\eta_i}\norm{\x_i^\star-\yt_i^1}^2} + \Psi^1_1\right]
	\end{align}
	If we start with $\beta_1 = 1/2$, then it can be shown that $\beta_s \leq \frac{2}{s+3}$ (see \cite{shang2018asvrg}), implying that
	\begin{align}
	\sum_{i\in\V}\Ex{D_{f_i+h_i}(\xt_i^{S+1},\x_i^\star)} &\leq \frac{2W_1}{(S+3)^2}\left[V + \frac{G_\rho}{m} + \sum_{i\in\V}\frac{1}{\eta_i m} \right]
	\end{align}
	Note that for this $\beta_s$, we also set $p_s = \beta_s$. Recall also that for $\omega_i^s = \frac{L_i}{\beta_s(1-\beta_s)}$, we require $\eta_i \leq \frac{1-\beta_s}{(1+\beta_s(1-\beta_s)) L_i}$. The bound on $\eta_i$ is the smallest for $s = 1$, so can be ensured by selecting $\eta_i  =  \frac{2}{5L_i}$ resulting in the bound
	\begin{align}
		\sum_{i\in\V}\Ex{D_{f_i+h_i}(\xt_i^{S+1},\x_i^\star)} &\leq \frac{2W_1}{(S+3)^2}\left[V + \frac{G_\rho}{m} + \frac{5}{2m}\sum_{i\in\V}L_i \right]
	\end{align}

	Let us simplify the bound by taking $L_i = L$ and $K_i = K$ for all $i\in \V$. Then, the number of samples required is $\O(SK+Sm)$ while the number of message exchanges required is $\mathcal{O}(S^2m)$.  Selecting $m = \O(K)$, we obtain the bound
	\begin{align}
		\sum_{i\in\V}\Ex{D_{f_i+h_i}(\xt_i^{S+1},\x_i^\star)} \leq \O\left(\frac{E+V(L+K)}{KS^2}\right)
	\end{align} 
	which translates to the oracle complexity
	\begin{align}
	\O\left(\frac{\sqrt{K(E+V(K+L))}}{\sqrt{\epsilon}}\right)
	\end{align} 
	and the communication complexity of $\O(\frac{E + V(K+L)}{\epsilon})$. 
\end{proof}}
As expected, the oracle complexity of the ASVR-PD algorithm is significantly better than that of the  non-accelerated variant.
\section{Numerical Tests}
This section provides numerical results and comparisons between the proposed and other state-of-the-art algorithms. We illustrate the performance of various algorithms on the binary classification problem. A regularized logistic regression objective is considered, where the loss function at node $i$ takes the form
\begin{align}
f_i(\x) = \frac{1}{K_i} \sum_{t=1}^{K_i} \log[1+ \exp(-c_{t}\d^T_{t}\x)] + \frac{\tau}{2}\norm{\x}^2\label{lr}
\end{align}
where $\tau > 0$ is the regularization parameter. The training samples comprise of feature vectors $\d_{t}\in \Rn^n$ and corresponding binary labels $c_{t}$ for $1\leq t \leq K_i$ and $i\in\V$.   


We make use of the Epsilon dataset \footnote{\url{https://www.csie.ntu.edu.tw/~cjlin/libsvmtools/datasets/binary.html}}, where $n = 2000$. Unless otherwise stated, we consider a fixed connected network with $V = 50$ nodes and $E = 250$ edges. Each node is assumed to have access to $K_i = 2000$ data points selected randomly from the dataset. 


We compare the performance of the proposed SVR-PD, SAGA-PD, LSVR-PD, and ASVR-PD algorithms with various state-of-the-art methods, namely, Stochastic Gradient Tracking (S-DIGing) \cite{li2019s}, Edge-based primal-dual algorithm (EDGE) \cite{wang2019edge}, Network distributed approximate Newton-type method (Network DANE) \cite{li2019communication}, Gradient Tracking Stochastic Variance Reduction (GT-SVR) \cite{xin2019variance}, Decentralized training over Decentralized data (D$^2$) \cite{tang2018d}, Exact First-Order Algorithm (EXTRA) \cite{shi2015extra}, and the Stochastic  Distributed Gradient Descent (DGD) \cite{lian2017can}. As the focus here is on the convergence rate, we plot the cumulative Bregman divergence defined in \eqref{bregman} against the number of calls to the stochastic gradient oracle. The optimal $\x^\star$ is obtained by running one of the algorithms (specifically SVR-PD) for several iterations till convergence.  We set the regularization parameter $\tau = 0.0014$ in order to obtain reasonable classification accuracy, though that is not the focus here.

Recall that the SVR-PD algorithm requires an exponentially increasing epoch size. As in \cite{allen2016improved}, we use the update rule $m_{s+1} = \min\{2m_s, \text{maxiter}\}$ where we have set maxiter $= 1000$. We also plot the performance of an early-stopping SVR-PD variant, where the inner iterations are terminated if the progress is below a certain threshold. Such an early-stopping SVR-PD was not analyzed but is nevertheless presented as a potential competing algorithm. In ASVR-PD, The starting values are set as $p_1 = 0.2$ and $\beta_1 = 0.99$. Other parameters for different algorithms were chosen manually and summarized in Table \ref{param_list}. The Algorithms GT-SVRG and Network Dane require doubly stochastic matrices that are generated using the metropolis strategy. The step-sizes for Network Dane is selected to be 10 and that for GT SVR is set to be 0.1, as these values yielded the fastest convergence rates for the respective algorithms. Along with these algorithms the proposed framework was also compared with  EDGE of \cite{wang2019edge} and the stochastic DGD from \cite{lian2017can}. The maximum number of iterations for all the algorithms is fixed to be 7000.

\begin{table}[h]	
\caption{Parameters for the proposed algorithms}\label{param_list}
	\centering
	{\renewcommand{\arraystretch}{1}
\begin{tabular}{|l|c|c|}
	\hline
	Algorithm & $\eta_i$ & $\rho$\\
	\hline
	Decentralized SAGA-PD & 0.5 & 0.5\\
	Decentralized SVR-PD & 0.7&0.9 \\
	Decentralized LSVR-PD & 0.6 & 0.5\\
	Decentralized Accelerated SVR-PD & 0.7 & 0.9\\
	\hline
\end{tabular}}
\end{table}

\begin{figure}[ht!]
\centering
\includegraphics[clip,trim={0.5cm 2cm 2cm 0.5cm},scale=0.63]{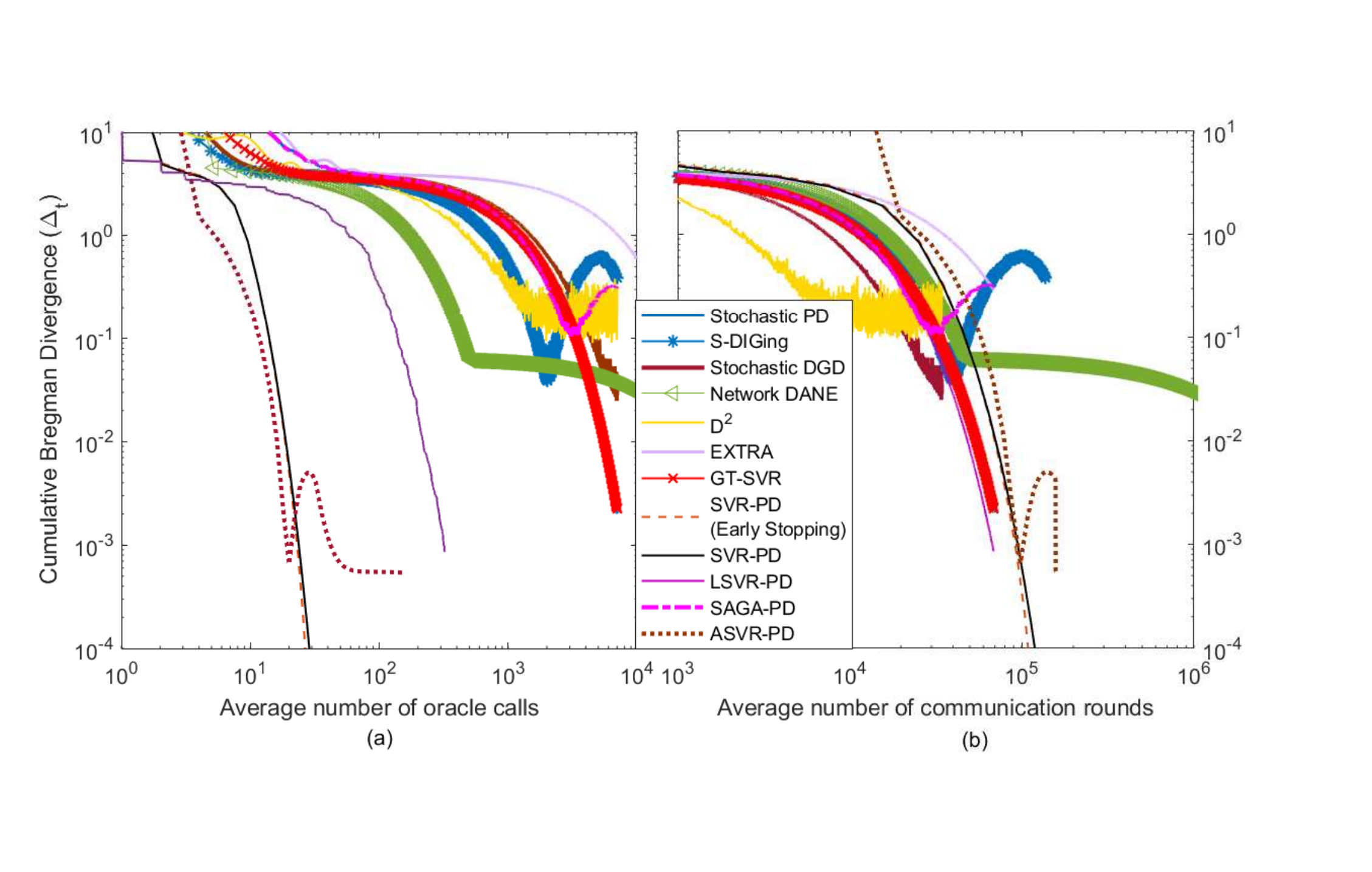}
\caption{Comparison of Oracle calls and communication rounds for proposed and state-of-the-art decentralized algorithms}	\label{breg_dive}
\end{figure}
Fig. \ref{breg_dive} shows the evolution of the Bregman divergence $\Delta_t$ against the average number of oracle calls as well as against the number of communication rounds. In terms of the oracle complexity, it can be seen that SVR-PD, SVR-PD with early stopping, and ASVR-PD are significantly faster than all other algorithms. Among the state-of-the-art algorithms, the Network DANE \cite{li2019communication} achieves the best performance but is far from some of the proposed algorithms. {We observed that for the EDGE \cite{wang2019edge} algorithm, we were unable to obtain a reasonable objective function value for the set of parameters tuned to maximize classification accuracy.}

Although the ASVR-PD is theoretically superior to all other algorithms, we found it difficult to tune in practice and it exhibited oscillations when the Bregman divergence was below a certain threshold. In terms of communication complexity, the performance of most of the algorithms is largely similar, though the proposed algorithm achieved lower values of $\Delta_t$ within the stipulated number of iterations. In summary, Fig. \ref{breg_dive} establish SVR-PD and SVR-PD with early stopping as the state-of-the-art algorithms for decentralized optimization. 

\begin{figure}[ht!]
	\centering
	\includegraphics[clip,trim={3cm 0cm 3.5cm 0cm},scale=0.5]{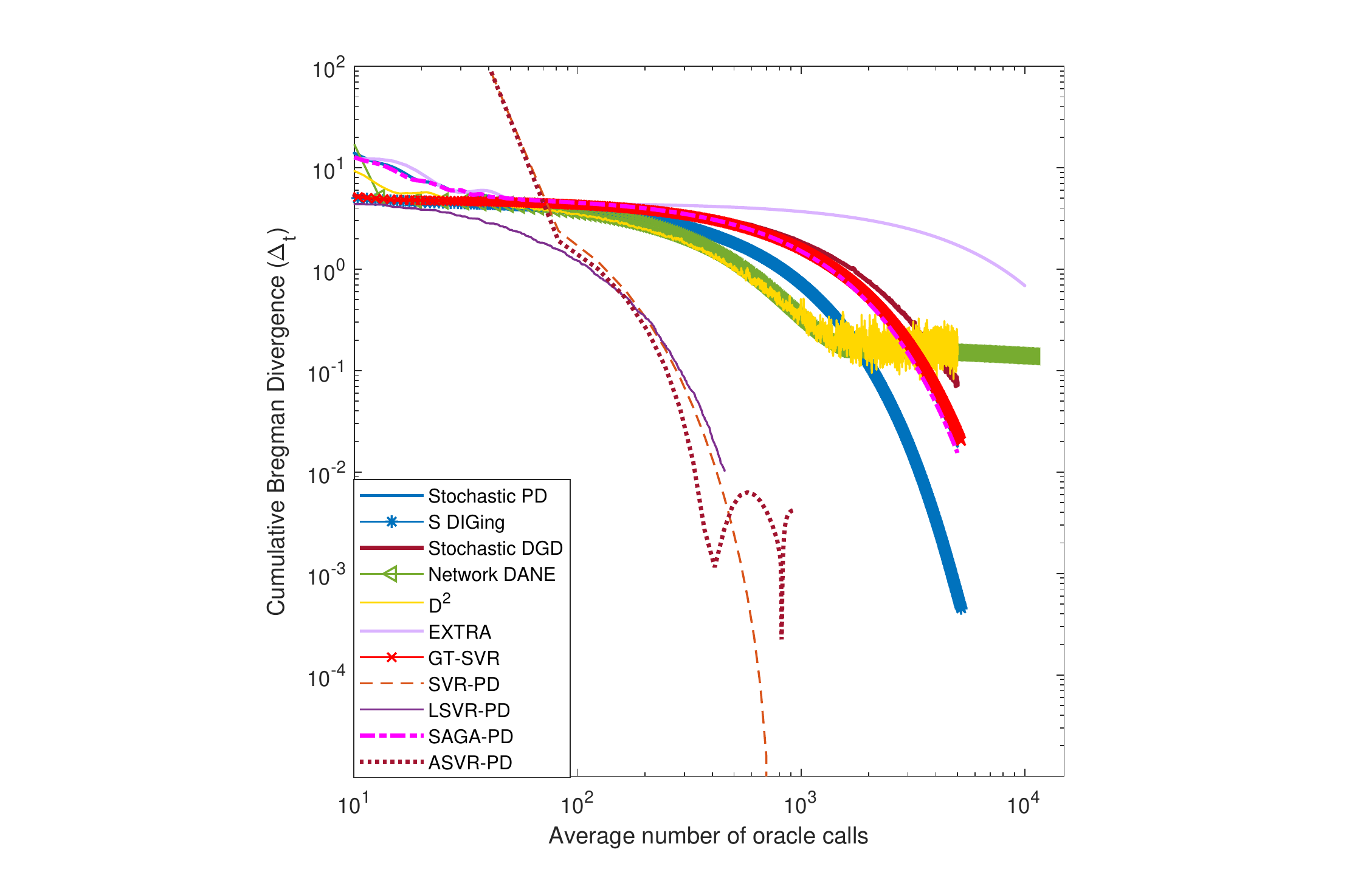}
	\caption{Performance comparison for $K_i = 50$.}
	\label{small_data}
\end{figure}

Recall that theoretically, the performance of some of the proposed algorithms did depend on the relative magnitudes of $K$ and $L$. To see if the value of $K$ has any actual impact on the performance, we consider the same setting but with only $50$ data points per node. The relative performance of most of the algorithm remains unchanged (Fig. \ref{small_data}), except for that of LSVR-PD whose performance improves significantly and becomes almost same as that of SVR-PD. 

Next, we study the impact of topology, regularization, and data-set size for various algorithms. Henceforth, we will only depict the performance of the SVR-PD with early stopping. It is remarked that for all the remaining figures, the performance of SVR-PD was found to be almost identical, while that of ASVR-PD was slightly better so long as $\Delta_t$ was approximately above $10^{-2}$. 

\subsection{Evaluation for large dataset}
Next, let us consider the opposite case, that is when the size of the dataset is much larger. We consider the same settings are earlier but with $K_i = 8000$ data points per node. For this case as well, SVR-PD continues to be the most efficient among all that are tested. Among the state of the art algorithms, the performance of GT-SVR and Network-DANE is the best, but worse than that of SVR-PD. 
\begin{figure}[ht!]
	\centering
	\includegraphics[clip,trim={3cm 0cm 3.5cm 0cm},scale=0.5]{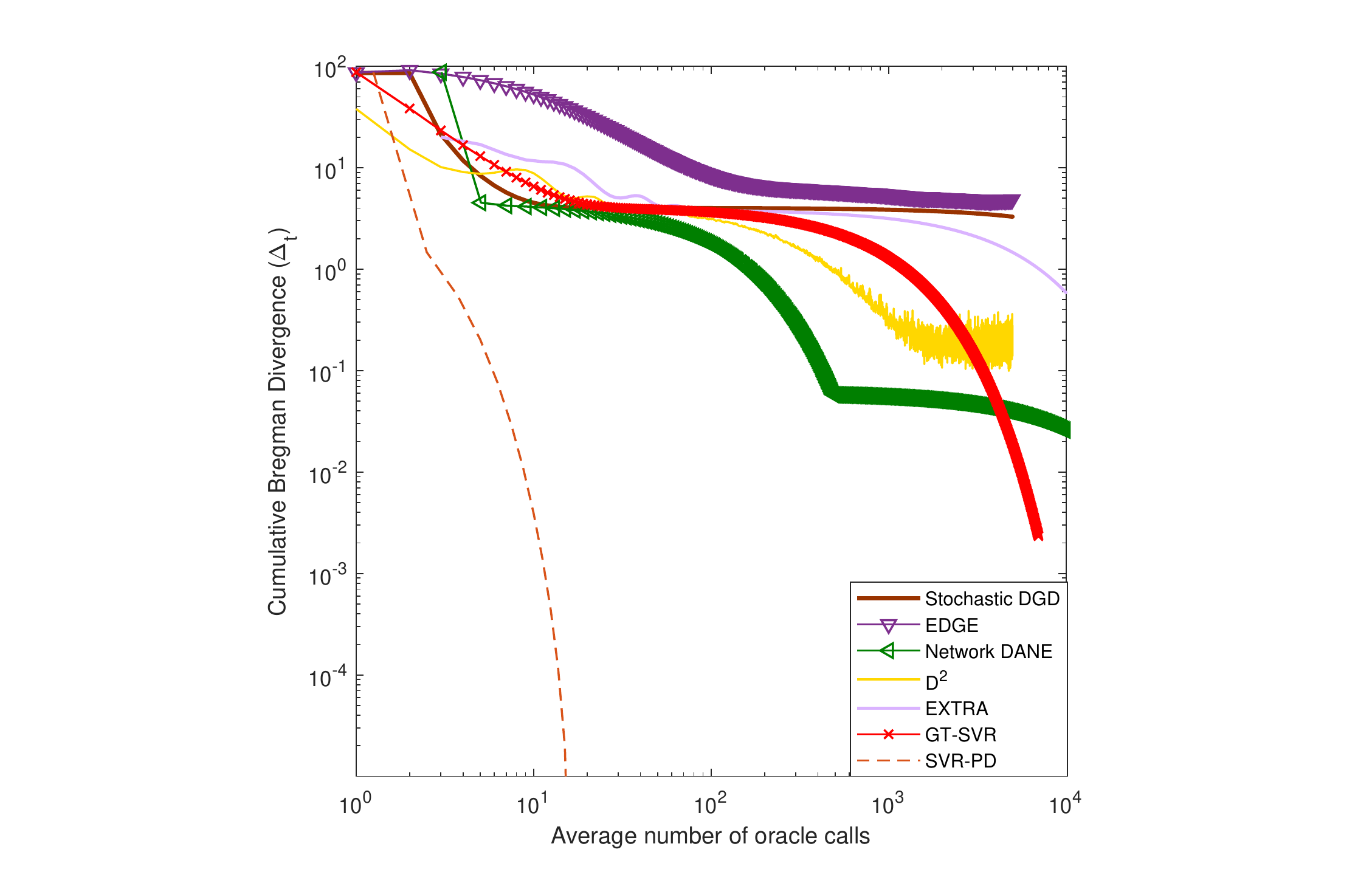}
	\caption{Performance comparison for $K_i = 8000$.}
	\label{bigdata}
\end{figure}

\subsection{Impact of topology}
\begin{figure}[ht!]
	\centering
	\includegraphics[clip,trim={3cm 0cm 2cm 0cm},scale=0.55]{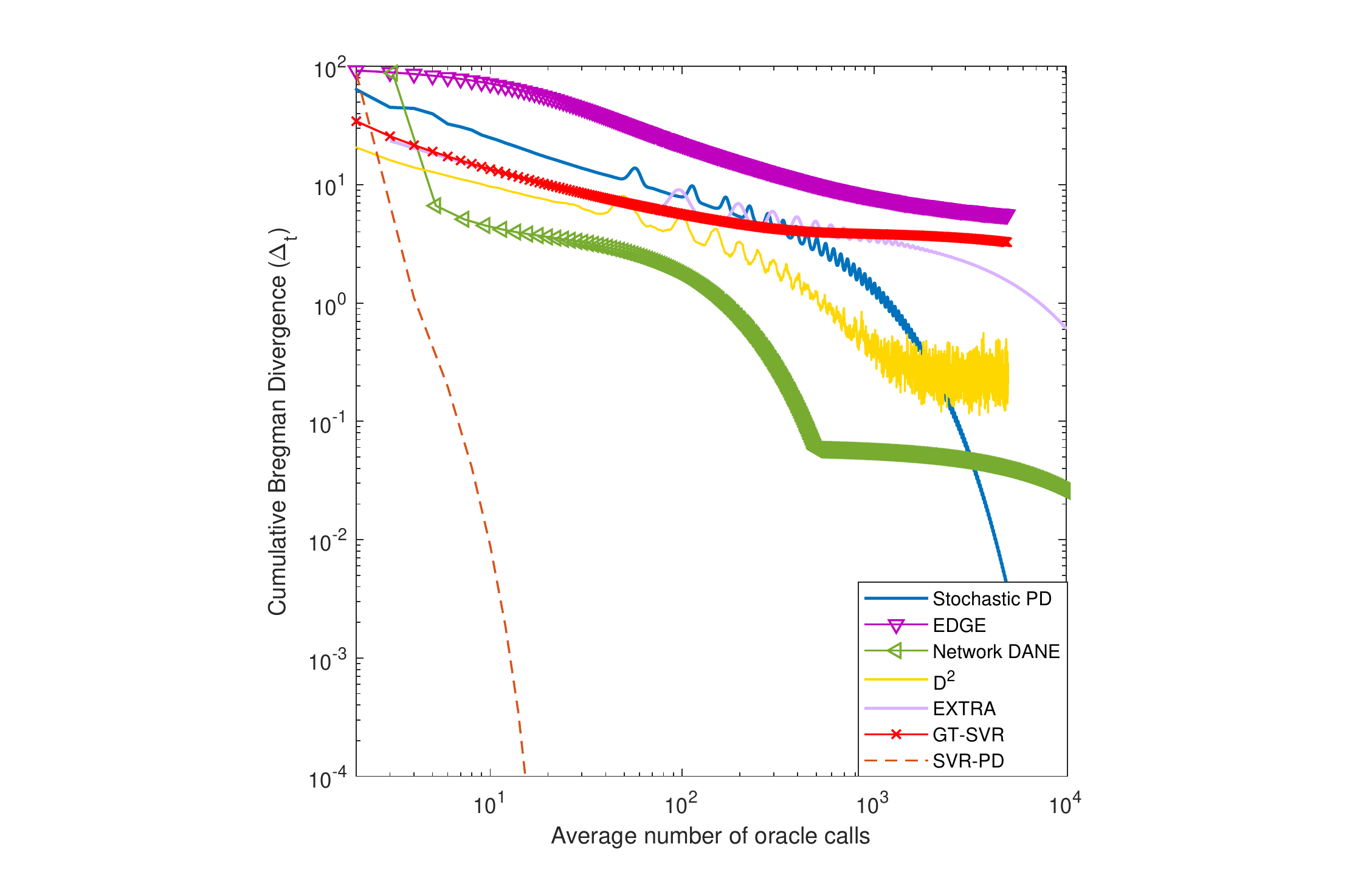}
	\caption{Performance comparison for ring topology}
	\label{ring}
\end{figure}

The performance of various algorithms is evaluated for the less well-connected ring topology with 50 nodes. As shown in Fig. \ref{ring}, the lack of good connectivity seems to have little impact on the performance of the SVR-PD algorithm, which is again the fastest and significantly better than the other proposed algorithms as well as the state-of-the-art algorithms. We do not include the performance of the S-DIGing and the Stochastic DGD algorithms, which do not converge for this topology. The value of $\rho$ is kept the same as that in Fig. 1.

\subsection{Logistic Regression without Regularization}
\begin{figure}[ht!]
	\centering
	\includegraphics[clip,trim={3cm 0cm 3.5cm 0cm},scale=0.5]{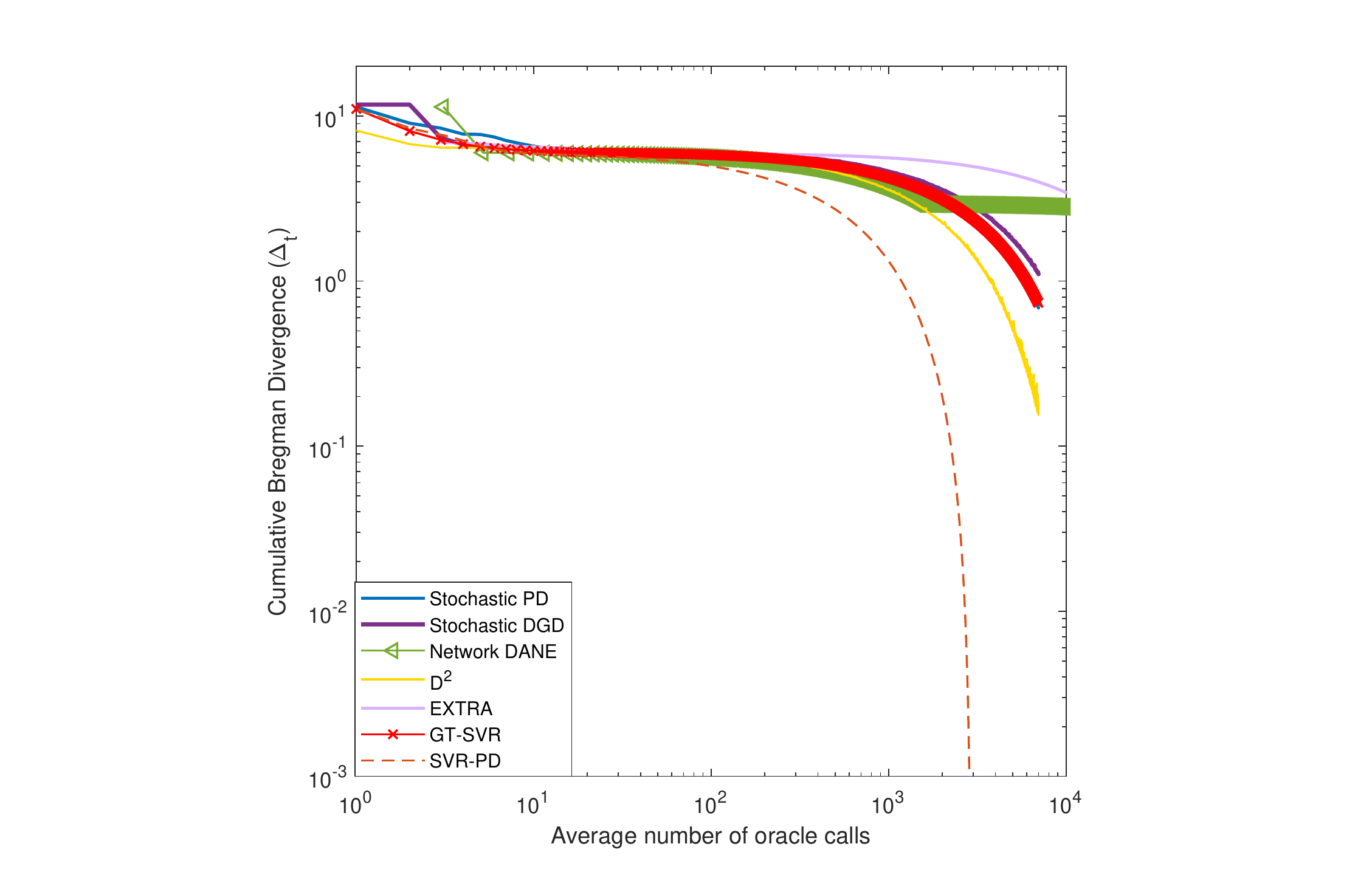}
	\caption{Performance comparison for classical (un-regularized) logistic loss function}
	\label{convex_fun}
\end{figure}

We examine the effect of not using the squared-norm regularizer in the logistic loss since the proposed class of algorithms is capable of handling non-strongly-convex objectives. Fig. \ref{convex_fun} shows the evolution of the Bregman divergence for the classical logistic loss function, while all other settings of Fig. \ref{breg_dive} are retained. As expected, the convergence rate for this case is much worse, and the performance of some algorithms suffers by more than two orders of magnitude. Nevertheless, the SVR-PD algorithm is still the fastest, while the $D^2$ algorithm is a distant second. The S-DIGing algorithm does not converge for this case, and its performance is not shown. 

\subsection{Run-time complexity}
\begin{figure}[ht!]
	\centering
	\includegraphics[scale=0.5]{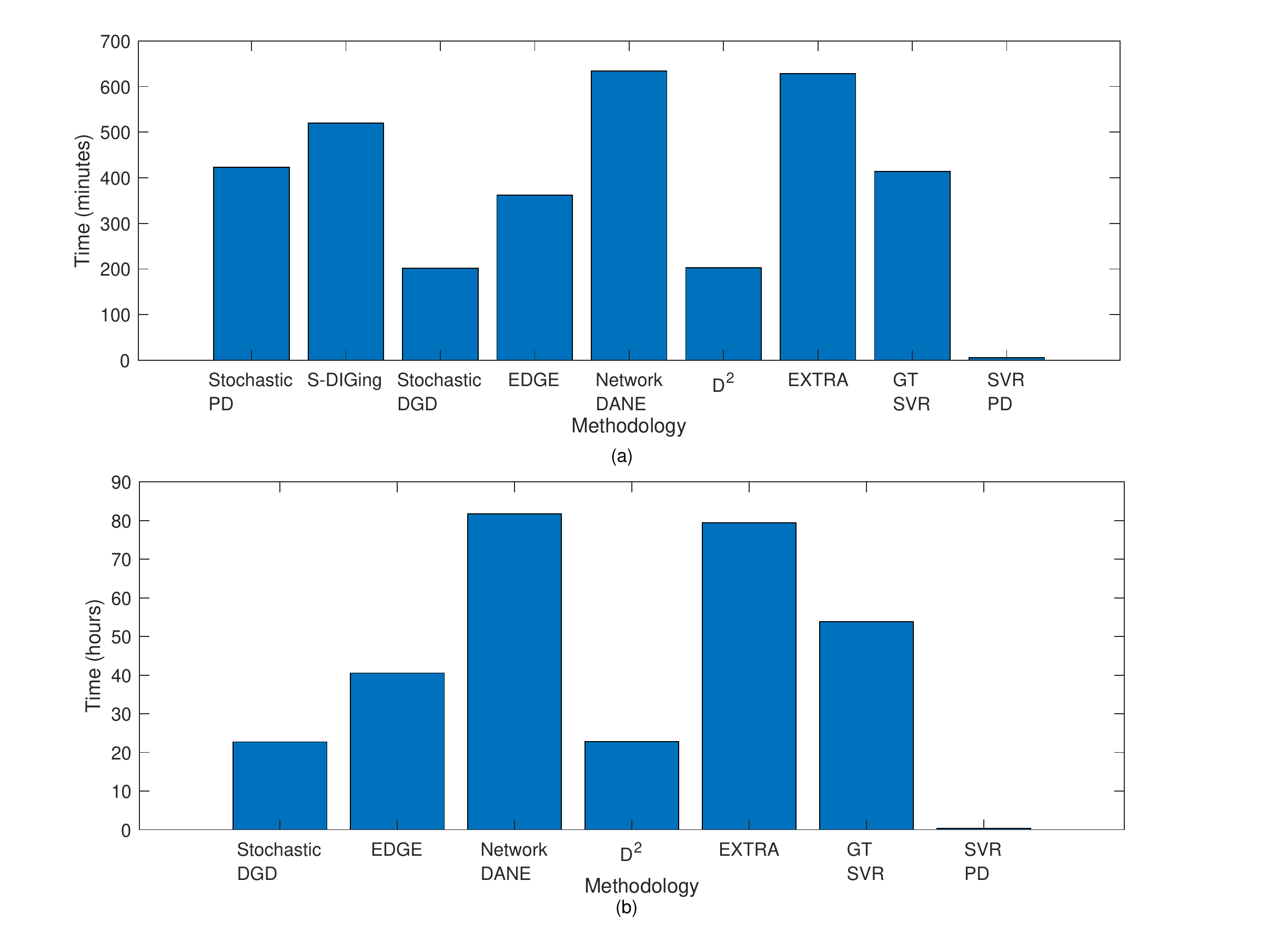}
	\caption{Run time comparison for (a) $K_i=2000$ (b) $K_i=8000$}
	\label{time}
\end{figure}

A more practical metric when comparing various algorithms can be their run-time or wall-clock time performance, though it is machine-dependent and therefore not as general as oracle complexity. To this end, we have recorded the run-time performance of various algorithms using a common implementation framework and on a single machine. Specifically, all the algorithms were implemented in MATLAB 2016a within a common wrapper function and all the tests were run on a machine with Intel(R) Celeron(R) M Processor with 64 GB of RAM. Fig. \ref{time}a shows the run-time of various algorithms for the parameters and settings of Fig. \ref{breg_dive}. Each algorithm is run till its optimality gap goes below $10^{-4}$ or the maximum number of iterations is reached. Here again, we observe the significant difference in the performance of the proposed SVR-PD algorithm and that of the state-of-the-art algorithms. As expected, SVR-PD is also the fastest in terms of the run-time. 

The same test was repeated for the case when $K_i = 8000$ and the run-time performance is shown in Fig. \ref{time}b. As earlier, the performance of SVR-PD remains the best, which yielded the result in approximately 30 minutes. In contrast, all other algorithms took 20-80 hours. 

\section{Conclusions}
\label{sec:conclusions}
This work proposed a primal-dual framework for the construction of decentralized stochastic algorithms. A generic decentralized stochastic primal-dual algorithm is proposed that can be adapted to all variance-reduction techniques such as SAGA, LSVRG, and SVRG++. Further, a decentralized accelerated version of the SVRG is proposed that is shown to achieve a complexity of $\O(\frac{1}{\sqrt{\epsilon}})$. Such results imply that the proposed class of decentralized algorithms is optimal since it achieves the centralized rate. 

Looking forward, it still remains to be seen if the proposed algorithms exhibit faster convergence under strong convexity. The proof approach developed here works only for unbiased stochastic gradient approximations, and extension to biased case also remains open. Finally, it would be interesting to see if the analysis developed here carries over to the non-convex case.
\appendix
\section{Proof of Lemma \ref{opt}}\label{exactgrad}
	Since $\X^\star$ is a solution to \eqref{P1}, it already satisfies \eqref{pfeas}. Let us re-write \eqref{P1} as 
	\begin{align}\label{P1a}
	\X^\star& = \arg\min_{\X} \sum_{i=1}^V [f_i(\x_i) + h_i(\x_i)] \\
	\text{s. t. }  & \x_i = \x_j \hspace{4mm} j\in\N_i, i\in\V \label{c1}\\
	& \x_j = \x_i	\hspace{4mm} j\in\N_i, i \in \V \label{c2}
	\end{align}
	where we have simply repeated each constraint. Associating dual variables $\bnu_{ij}$ and $\bnu_{ji}$ with the constraints in \eqref{c1} and \eqref{c2}, respectively, the KKT conditions of \eqref{P1a} are given by 
	\begin{align}
	\nabla f_i(\x_i^\star) + \v_i^\star +  \sum_{j\in\N_i} (\bnu_{ij}^\star - \bnu_{ji}^\star) &= 0 \label{foonu}
	\end{align}
	and $\x_i^\star = \x_j^\star$. Observe that \eqref{foonu} is of the required form in \eqref{fo} for the choice $\lam_{ji}^\star = \bnu_{ij}^\star - \bnu_{ji}^\star$ for all $j\in\N_i$ and $i\in \V$. For this choice, it also holds that $\lam_{ij}^\star = \bnu_{ji}^\star - \bnu_{ij}^\star = -\lam_{ji}^\star$ which is the required condition in \eqref{dfeas}.

	Next, we establish that $(\X^\star,\Lam^\star)$ is a saddle point of the augmented Lagrangian $L(\X,\Lam)$. Note that $L(\X,\Lam)$ is a affine function of $\Lam$ defined over $\Lambda$ so that for $\max_{\Lam \in \Lambda}L(\X^\star,\Lam)$ to be bounded, it is required that
	$\x_i^\star = \x_j^\star $ for all $j\in\N_i$ and $i\in \V$. In this case, it can be seen that $L(\X^\star,\Lam) = L(\X^\star,\Lam^\star)$ for any $\Lam\in\Lambda$. Further, the condition $\X^\star = \min_{\X}L(\X,\Lam^\star)$ is equivalent to
	\begin{align}
	\nabla f_i(\x_i^\star) + \v_i^\star + \sum_{j\in\N_i}\lam_{ji}^\star + \rho\sum_{j\in\N_i} (\x_i^\star - \x_j^\star) &= 0
	\end{align}
	which for $\x_i^\star = \x_j^\star$ is exactly the same as \eqref{fo} and thus holds. Hence we conclude that $(\X^\star,\Lam^\star)$ with $\Lam^\star$ defined as earlier is a saddle point of the augmented Lagrangian $L(\X,\Lam)$. 

Let $(\Xt,\Lt)$ be limit points of the iterations \eqref{primup1}-\eqref{dualup1} in the exact gradient case, i.e., when $\g_i^t = \nabla f_i(\x_i^t)$. Then these points satisfy:
\begin{align}
\vt_i + \nabla f_i(\xt_i) + \sum_{j\in\N_i}\lt_{ji} = 0 \label{lim1}\\
\lt_{ij}+\lt_{ji} - \rho(\xt_i - \xt_j) = 0\label{lim2} 
\end{align}
where $\vt_i \in \partial h_i(\xt_i)$. Note that \eqref{lim2} is written in two different ways for every $(i,j)\in\E$:
\begin{align}
\lt_{ij}+\lt_{ji} - \rho(\xt_i - \xt_j) &= 0 \\
\lt_{ji}+\lt_{ij} - \rho(\xt_j - \xt_i) &= 0 
\end{align}
corresponding to updates at node $i$ and $j$ respectively. Adding, we obtain $\lt_{ij} + \lt_{ji} = 0$, which implies that $\xt_i = \xt_j$ since $\rho > 0$. Since the graph $\G$ is connected, let us denote $\xt = \xt_i$ for any $i\in\V$. Summing \eqref{lim1} over all $i \in \V$, and substituting, we obtain
\begin{align}
-\sum_{i\in\V}\nabla f_i(\xt) \in \partial(\sum_{i\in\V} h_i(\xt))
\end{align}
since $\sum_{i\in\V}\sum_{j\in\N_i} \lt_{ji} = 0$. In other words, the limit point $\xt$ is the solution to \eqref{P1}.

\section{Link with PDMM Algorithm}\label{PDMM}
The PDMM algorithm was first proposed in \cite{main} and considered the affine-constrained deterministic optimization problem
\begin{align}\label{pdmmprob}
\X^\star& = \arg\min_{\X} \sum_{i\in \V} f_i(\x_i) \\
& \text{s. t. }  \  \A_{ij}\x_i+\A_{ji}\x_j=\c_{ij} \hspace{4mm} \forall(i,j) \in \E \nonumber.
\end{align}
where $\A_{ij} \in \Rn^{K\times n}$ and $\c_{ij} \in \Rn^K$. As a special case, consider $\c_{ij} = 0$, $\A_{ij}  = -\A_{ji} = a_{ij}\mathbf{I}$ for some $a_{ij} \in \{-1,1\}$, in which case, \eqref{pdmmprob} becomes a special case of \eqref{P1} for $h_i = 0$. The general PDMM algorithm also includes a parameter matrix $\P_{ij}$ which we set equal to $\rho\I$. With these specializations, the quadratic approximation PDMM algorithm proposed in \cite{qapdmm} entails carrying out the updates:
\begin{align}
\x_i^{t+1} &= \left(\tfrac{1}{\eta_i} + \rho\abs{\N_i} \right)^{-1}\left[\tfrac{1}{\eta_i}\x_i^t - \nabla f_i(\x_i^t) + \sum_{j\in\N_i}a_{ij}\bar{\lam}_{ji}^t + \rho\x_j^t\right] \label{pdmm-xup}\\
\bar{\lam}_{ij}^{t+1} &= \bar{\lam}_{ji}^t + \rho a_{ij}(\x_j^t-\x_i^{t+1}) \label{pdmm-lup}
\end{align}
Multiplying \eqref{pdmm-lup} by $a_{ij}$ and substituting $\lam_{ij}^t = a_{ij}\bar{\lam}_{ij}^t$ and $\lam_{ji}^t = -a_{ij}\bar{\lam}_{ji}^t$, we obtain the proposed updates \eqref{primup}-\eqref{dualup}:
\begin{align}
\x_i^{t+1} &= \left(\tfrac{1}{\eta_i} + \rho\abs{\N_i} \right)^{-1}\left[\tfrac{1}{\eta_i}\x_i^t - \nabla f_i(\x_i^t) - \sum_{j\in\N_i}\lam_{ji}^t + \rho\x_j^t\right] \label{pdmm-xupf}\\
\lam_{ij}^{t+1} &= -\lam_{ji}^t + \rho (\x_j^t-\x_i^{t+1}). \label{pdmm-lupf}
\end{align}
The connection suggests that the proposed algorithms can be generalized to allow general affine constraints such as those in \eqref{pdmmprob}. Conversely, quadratic approximation PDMM can be generalized to solve non-smooth composite problems.

\section{Decentralized Sketched Primal-Dual}\label{pdmmsec}
We demonstrate the applicability of the proposed algorithms to settings beyond variance-reduced methods. Consider the high-dimensional setting, i.e., where $n$ is large. As in the earlier sections, the high-dimensional setting precludes the possibility of calculating the full gradient at every iteration. Instead, each call to an oracle returns a random component $[\nabla f_i(\x)]_j$ for some $1\leq j\leq n$. Such a gradient component is, for instance, returned by a zeroth-order oracle, where 
\begin{align}
[\nabla f_i(\x)]_j \approx \frac{1}{\varepsilon}(f(\x+\varepsilon \e_j) - f(\x))
\end{align}
for a small $\varepsilon > 0$,  where $\e_j$ is the vector with one at the $j$-th location and zeros elsewhere. Zeroth order stochastic gradient approaches have been widely studied. Here, we consider a variant of the SEGA algorithm considered in \cite{gorbunov2019unified}. The proposed coordinate-based PD algorithm is given in Algorithm \ref{algosega}. Here $\odot$ denotes the entrywise product and the expressions for $\h_i^{t+1}$ and $\g_i^t$ can equivalently be written as
\begin{align}
[\h_i^{t+1}]_j &= \begin{cases} [\nabla f_i(\x_i^t)]_j	& j = j_i^t \\
[\h_i^t]_j	& j\neq j_i^t
\end{cases}\\
[\g_i^t]_j &= \begin{cases}
n[\nabla f_i(\x_i^t)]_j+ (1-n)[\h_i^t]_j & j = j_i^t \\
[\h_i^t]_j	& j \neq i_t.
\end{cases}
\end{align}
Observe that the algorithm maintains two $n \times 1$ vectors $\h_i$ and $\g_i$ at each node $i \in\V$. However, only a single component of these vectors is updated at each iteration. 
\begin{algorithm}
	\caption{SkPD}\label{algosega}
	\begin{algorithmic}
		\STATE Initialize $\{\x_i^0\}$,  $\{\lam^0_{ij}\}$, and $\{\h_i^0 = 0\}$ 
		\FOR { $\textit{t} = 1,2,... $ } 
		\STATE pick $j_i^t \in \{1, \ldots, n\}$ uniformly at random and obtain $[\nabla f_i(\x_i^t)]_{j_i^t}$ for all $i\in\V$,
		\STATE update $\h_i^{t+1} = \h_i^t + \e_{j_i^t} \odot(\nabla f_i(\x_i^t) - \h_i^t)$ for all $i\in\V$,
		\STATE set $\g_i^t = \h_i^t + n\e_{j_i^t}\odot(\nabla f_i(\x_i^t) - \h_i^t)$ for all $i\in\V$,
		\STATE update $\x_i^{t+1}$ using \eqref{primup} for all $i\in\V$,
		\STATE update $\lam_{ij}^{t+1}$ using \eqref{dualup} for all $i\in\V$, $j\in\N_i$	 		
		\ENDFOR
	\end{algorithmic}
\end{algorithm}

The following lemma provides the oracle complexity result for the SkPD algorithm. Here, a sketched gradient oracle call returns a random component of $\nabla f_i(\x)$ for a given $\x$. 

\begin{lemma}\label{sega-pd}
	The SkPD algorithm incurs a sketched-gradient oracle complexity of $\O\left(\frac{E+V(1+nL)}{\epsilon}\right)$ for $\eta = \frac{1}{8nL}$ and $L \geq \frac{1}{4\sqrt{2}}$. 
\end{lemma}
\begin{proof}
	The first part of the proof follows exactly as in \cite{gorbunov2019unified}[Lemma A.8], giving $A_i = 2nL_i$, $B_i = 2n$, $\varrho_i = 1/n$, and $C_i = L_i/n$ while requiring that $\eta_i \leq \frac{1}{8nL_i}$. From Lemma \ref{lema2}, we obtain
	\begin{align}
	\frac{1}{T}\sum_{t=1}^T \Delta^t & \leq \frac{W_1G_\rho}{T} + \frac{W_1}{T}\sum_{i\in \V}\left[\frac{1}{\eta_i} + 2n^2\eta_i\right]
	\end{align}
	To simplify this bound, let us set $L_i = L$ for all $i$ and assume that $L \geq \frac{1}{4\sqrt{2}}$. Then for $\eta = \frac{1}{8nL}$, we obtain the bound
	\begin{align}
		\min_{1\leq t \leq T}\Delta^t &\leq \O\left(\frac{E+V(1+nL)}{T}\right)
	\end{align}
	which translates to the required oracle complexity. 
\end{proof}

\bibliographystyle{siamplain}
\bibliography{ref}
\end{document}